\documentclass[12pt]{amsart}
\usepackage{geometry}
\usepackage{amsmath,amsfonts,amssymb,graphicx,hyperref, amsthm,mathabx,amsrefs,geometry}
\hypersetup{colorlinks=true,citecolor=blue,linkcolor=blue,urlcolor=blue,pdfstartview=FitH,bookmarksdepth=3}
\usepackage[noabbrev, nameinlink]{cleveref}
\usepackage{enumitem}
\setlist[enumerate]{label=\arabic*.}

\newcommand{\paren}[1]{\ensuremath{\left( #1 \right)}}
\newcommand{\set}[1]{\ensuremath{\left\{ #1 \right\}}}
\newcommand{\abs}[1]{\ensuremath{\left| #1 \right|}}
\newcommand{\norm}[1]{\ensuremath{\left\| #1 \right\|}}
\newcommand{\setdiv}{\,\middle|\,}
\newcommand{\summod}[1]{\ensuremath{\,(\mathrm{mod}\,#1)}}

\newcommand{\Matrix}[1]{\begin{pmatrix}#1\end{pmatrix}}

\newcommand{\piecewise}[1]{\left\{\begin{matrix}#1\end{matrix}\right.}
\newcommand{\If}{\mbox{if }}

\renewcommand{\Re}{{\mathop{\mathgroup\symoperators Re}}}
\renewcommand{\Im}{{\mathop{\mathgroup\symoperators Im}}}
\newcommand{\sgn}{{\mathop{\mathgroup\symoperators \,sgn}}}

\newcommand{\Max}[1]{\ensuremath{\max \set{#1}}}
\newcommand{\Min}[1]{\ensuremath{\min \set{#1}}}
\newcommand{\floor}[1]{\ensuremath{\left\lfloor #1 \right\rfloor}}
\newcommand{\ceil}[1]{\ensuremath{\left\lceil #1 \right\rceil}}

\newcommand{\Z}{\mathbb{Z}}
\newcommand{\R}{\mathbb{R}}
\newcommand{\N}{\mathbb{N}}

\newcommand{\C}{\mathbb{C}}

\newcommand{\wbar}[1]{\overline{#1}}
\newcommand{\wtilde}[1]{\widetilde{#1}}
\newcommand{\what}[1]{\widehat{#1}}

\newcommand{\BigO}[2][]{O_{#1}\paren{#2}}
\newcommand{\e}[1]{e\paren{#1}}
\newcommand{\trans}[1]{{#1}^T}

\DeclareMathOperator{\Tr}{Tr}
\DeclareMathOperator*{\res}{res}

\theoremstyle{plain}
\newtheorem{thm}{Theorem}

\newtheorem{lem}[thm]{Lemma}
\newtheorem{prop}[thm]{Proposition}

\crefname{thm}{Theorem}{Theorems}
\crefname{lem}{Lemma}{Lemmas}
\crefname{prop}{Proposition}{Propositions}
\crefname{cor}{Corollary}{Corollaries}
\crefname{conj}{Conjecture}{Conjectures}

\newcommand{\AdSq}{{\mathrm{Ad}}^2}

\newcommand{\sinmu}{\mathbf{sin}}
\newcommand{\cscmu}{\mathbf{csc}}

\DeclareMathOperator{\arcsinh}{arcsinh}
\newcommand{\Weyl}{W}

\title{On sums of hyper-Kloosterman Sums}
\author{Jack Buttcane}
\date{27 May 2022}

\begin{document}

\begin{abstract}
A formula of Kuznetsov allows one to interpret a smooth sum of Kloosterman sums as a sum over the spectrum of $GL(2)$ automorphic forms.
In this paper, we construct a similar formula for the first hyper-Kloosterman sums using $GL(3)$ automorphic forms, resolving a long-standing problem of Bump, Friedberg and Goldfeld.
Along the way, we develop what are apparently new bounds for the order derivatives of the classical $J$-Bessel function, and we conclude with a discussion of the original method of Bump, Friedberg and Goldfeld.
\end{abstract}


\maketitle

In 1988, Bump, Friedberg and Goldfeld (BFG) published a paper \cite{BFG} detailing the Fourier and spectral expansions of a spherical Poincar\'e series for $SL(3,\Z)$.
At the end of that paper, they theorized that one could isolate the contribution of the hyper-Kloosterman sum terms from that of the remaining Weyl cells.
Their particular method does not appear to be effective, but this idea has led to the long-standing goal of providing a spectral interpretation for a sum of hyper-Kloosterman sums.
In this paper, we provide one possible spectral expansion; the representation is not unique.

The desire to study such Poincar\'e series comes from the success of Kuznetsov's formulas for $SL(2,\Z)$.
Kuznetsov was attempting to address a conjecture of Linnik \cite{L01} on the average of the classical Kloosterman sums
\[ \operatorname{Kl}_2(m,n,c) = \sum_{\substack{x,y \summod{c}\\ xy\equiv 1\summod{c}}} \e{\frac{m x+n y}{c}}, \qquad \e{x} := \exp(2\pi i x), \]
which states that there should be good cancellation among the terms of an average over the moduli $c$.
Kuznetsov solved this problem in \cite{Kuz01} by relating the spectral and Bruhat expansions of (Fourier coefficients of) Poincar\'e series over $SL(2,\Z)$ to obtain a pair of formulas of the shape
\[ \sum_{c=1}^\infty \frac{\operatorname{Kl}_2(m,n,c)}{c} F_1\paren{\frac{\sqrt{mn}}{c}}+\text{trivial term} = \sum_\phi \frac{\lambda_\phi(n)\wbar{\lambda_\phi(m)}}{L(1,\AdSq \phi)} F_2(\mu_\phi)+\text{continuous term}, \]
where the sum on the right is over a basis of (minimal-weight) Maass cusp forms for $SL(2,\Z)$ with Hecke eigenvalues $\lambda_\phi(m)$ and spectral parameters $\mu_\phi$, and the (mostly) arbitrary test functions $F_1, F_2$ are related by Bessel or inverse-Bessel transforms (giving two formulas).
Kuznetsov's formulas are now fundamental to modern number theory as one of very few examples where we can detect cancellation in an average over the moduli of an exponential sum.

The paper \cite{BFG} of BFG determined explicitly the Bruhat decomposition of (Fourier coefficients of) an $SL(3,\Z)$ Poincar\'e series, which now contains not just one, but several different exponential sums corresponding to the elements of the Weyl group (essentially the permutation group on three letters).
The papers \cite{SpectralKuz, WeylI, WeylII} used this to construct the generalizations of the spectral (i.e. arbitrary test functions on the spectral side) Kuznetsov formula for $SL(3,\Z)$.

At the end of \cite{BFG}, the authors conjectured that a properly constructed kernel function would isolate the different Weyl-element terms in the Fourier coefficients of the corresponding Poincar\'e series; we address their proposed method in \cref{sect:BFGmethod}.
The paper \cite{ArithKuzII} was able to accomplish this for the long Weyl element, using an idea of Cogdell and Piatetski-Shapiro \cite{CPS}, the end result being a Kuznetsov-type formula having the (reasonably) arbitrary test function on the sum of long-element Kloosterman sums and no other Weyl element terms.
(The paper \cite{ArithKuzI} did this first, in a special case and using a different tactic.)

The long-element arithmetic Kuznetsov formula of \cite{ArithKuzII} is a significant advancement in that it shows we can construct such formulas outside the setting of rank-one groups (i.e. $SL(2)$ vs a rank-two group such as $SL(3)$), but the sheer complexity of the long-element Kloosterman sum makes it difficult to believe one could find such sums ``in the wild'' (i.e. outside the context of $SL(3,\Z)$ Poincar\'e series).
The more arithmetically interesting sums correspond to the order-three Weyl elements $w_4$ and $w_5$.
The $w_5$ Kloosterman sum (in the notation of \cite{ArithKuzII}) is
\begin{align*}
	S_{w_5}(m,n,c) =& \delta_{\substack{m_1 c_2=n_2 c_1^2 \\ c_1|c_2}} \sum_{\substack{x_1 \summod{c_1}, x_2 \summod{c_2}\\(x_1,c_1)=(x_2,c_2/c_1)=1}} \e{m_1\frac{\bar{x_1}x_2}{c_1}+m_2\frac{\bar{x_2}}{c_2/c_1}+n_1\frac{x_1}{c_1}},
\end{align*}
for $m,n,c\in\Z^2$, $c_1 c_2 \ne 0$.
In particular, at $n_2=m_1$,
\[ S_{w_5}(m,n,c) = c_1 \, \delta_{c_2=c_1^2} \sum_{\substack{x,y,z \summod{c_1}\\xyz\equiv1\summod{c_1}}} \e{\frac{n_1 x+m_1 y+m_2 z}{c_1}}, \]
and the three-variable sum at right is $\operatorname{Kl}_3(n_1,m_1,m_2,c)$, the first hyper-Kloosterman sum.

The hyper-Kloosterman sums $\operatorname{Kl}_3, \operatorname{Kl}_4,\ldots$ are relatively easy to locate in the context of other problems, and these have been studied, e.g. in \cite{Friedberg01, HNY01}.
The $w_4$ Kloosterman sums are symmetric to the $w_5$ sums, and the Kloosterman sums for the remaining Weyl elements are degenerate (i.e. either trivial or $\operatorname{Kl}_2$ sums).

In this paper, we achieve a spectral interpretation for a sum of the $w_5$ Kloosterman sums (hence also for a sum of $\operatorname{Kl}_3$ sums) over the moduli:
Writing
\[ \sum_{d=0}^\infty \int_{\mathcal{B}_3^{d*}} \ldots d_H\Xi \]
as shorthand for the full spectral expansion in $L^2(SL(3,\Z)\backslash SL(3,\R))$ (excluding the constant function and residual spectrum) over bases of minimal-weight forms (see \cref{sect:FormalPoin}), we have
\begin{thm}
\label{thm:MainCor}
Suppose $f:\R^\times \to \C$ is smooth and compactly supported outside the interval $(-T_1,T_1)$ and $m, n \in \Z_{\ne 0}^2$, then
\begin{align}
\label{eq:MainCor}
	\sum_{\substack{0 \ne c_1, c_2 \in \Z \\ m_1 c_2=n_2 c_1^2}} S_{w_5}(m,n,c) f\paren{C^3\frac{m_1^2 m_2 n_1}{n_2 c_1^3}} = \abs{m_1 m_2 n_1} \sum_{d=0}^\infty \int_{\mathcal{B}_3^{d*}} \frac{\lambda_\Xi(n)\wbar{\lambda_\Xi(m)}}{L(1,\AdSq\Xi)} \what{F}^d(\mu_\Xi) d_H\Xi,
\end{align}
where $\what{F}^d(\mu) = \what{F}^d(\mu, f, T_1)$ is made explicit in the next section.
\end{thm}

This formula differs from the previous such formulas in several ways.
First and foremost, the spectral interpretation is not unique; the function $\what{F}^d(\mu)$ depends on the choice of two auxiliary functions $h_1$ and $h_2$ and a parameter $X$ that must be taken large relative to $C T_1^{-1/3}$.
The spectral sum concentrates near the self-dual cusp forms, i.e. the symmetric squares of $SL(2,\Z)$ cusp forms.
Specifically, the main contribution comes from cusp forms whose spectral parameters are (up to permutation) of the form $(a+r,-a+r,-2r)$ with $a \asymp X^{\frac{1}{2}}$ and $r \ll X^\epsilon$; note that the self-dual forms occur at exactly $r=0$.
(The parameter $X$ does not appear on the left-hand side of \eqref{eq:MainCor}, but does affect the shape of the sum on the right-hand side; this is part of what we mean by ``the spectral interpretation is not unique''.)
This behavior was predicted in \cite{ArithKuzII} (see problems I, III and IV there).

From the $SL(3,\Z)$ Weyl law, the number of such forms is $X^{2+\epsilon}$ and the function $\what{F}^d(\mu)$ is of size $C^3 X^{-3/2}$ on this region, so the trivial bound on the spectral side has become $C^3 X^{\frac{1}{2}+\epsilon} \gg C^{\frac{7}{2}+\epsilon}$.
Meanwhile, the Deligne bound on the $w_5$ Kloosterman sum is roughly $c^{2+\epsilon}$, so the trivial bound on the Kloosterman sum side is $C^{3+\epsilon}$.
This is the second difference from previous formulas, as the spectral sums there were typically of length $C^\epsilon$.
This increase in the trivial bound will certainly increase the difficulty in applying the formula \cref{thm:MainCor} (or \cref{thm:MainThm}), but we expect applications still exist as it is an exact identity and we can apply the usual tools of analytic number theory (averaging, Poisson summation, etc.) to it; we give more details in the next section.

Finally, a related question is whether an arbitrary ``nice'' function can be expanded into Bessel functions.
We answered this question for the long-element Bessel function in \cite{ArithKuzII}, but we can extend this now to the Bessel functions attached to the other Weyl elements.
In particular, we have
\begin{thm}
\label{thm:w5BesselExpand}
For $f:\R^\times \to \C$ smooth and compactly supported,
\[ f(y) = \sum_{d \ge 0} \int_{\frak{a}^d_0} \what{F}^d(\mu) \, K^d_{w_5}((1,y),\mu) \sinmu^{d*}(\mu) d\mu \]
where $\what{F}^d(\mu) = \what{F}^d(\mu, f)$ has super-polynomial decay in $\mu$.
\end{thm}
Thus the $K^d_{w_5}(y,\mu)$ Bessel functions give a basis for the smooth, compactly supported functions on $\R^\times$.
(Though the notation suggests otherwise, the $K^d_{w_5}(y,\mu)$ Bessel function is only properly defined at $y_1=1$.)
However, it seems unlikely that the coefficients $\what{F}^d(\mu)$ are uniquely determined by the inner product
\[ \int_{-\infty}^\infty f(y_2) \wbar{K^d_{w_5}((1,y_2),\mu)} \frac{dy_2}{\abs{y_2}^3}, \]
as they are for the long-element Bessel function.

Expansions using the Bessel functions attached to the remaining Weyl elements are now known either classically or by symmetry with $K^d_{w_5}(y,\mu)$.
The proof of \cref{thm:w5BesselExpand} is an immediate application of the more general \cref{prop:Density}, below, to the test function constructed in \cref{sect:TestFun}.

\section{Results}
We require two auxiliary functions $h_1,h_2:\R \to \C$ which are to be smooth with $h_1$ supported on $(\frac{1}{2},2)$ and $h_1(t_1)=1$ for $t_1 \in (1-\frac{1}{10},1+\frac{1}{10})$; $h_2$ even, compactly supported on $(-\frac{1}{200},\frac{1}{200})$ and $h_2(0)=1$; these will affect the shape of the spectral terms and long-element Kloosterman sum terms, but crucially not the hyper-Kloosterman sum term.
Then the main theorem of this paper is
\begin{thm}
\label{thm:MainThm}
For $f,h_1,h_2$ as above,
\begin{align*}
	\abs{n_2} \sum_{\substack{0 \ne c_1, c_2 \in \Z \\ m_1 c_2=n_2 c_1^2}} S_{w_5}(m,n,c) f\paren{\frac{m_1^2 m_2 n_1}{n_2 c_1^3}} =& \abs{m_1 m_2 n_1 n_2} \sum_{d=0}^\infty \int_{\mathcal{B}_3^{d*}} \frac{\lambda_\Xi(n)\wbar{\lambda_\Xi(m)}}{L(1,\AdSq\Xi)} \what{F}^d(\mu_\Xi) d_H\Xi \\
	& -\sum_{0 \ne c_1, c_2 \in \Z} S_{w_l}(m,n,c) T_{w_l}(F)\paren{\frac{m_1 n_2 c_2}{c_1^2}, \frac{m_2 n_1 c_1}{c_2^2}},
\end{align*}
where
\begin{align*}
	T_{w_l}(F)(y) =& \delta_{y_1<0} \sqrt{\abs{y_1}} \sum_{\varepsilon \in \set{\pm1}} f(\varepsilon y_2 \sqrt{\abs{y_1}}) \wtilde{h_1}(\varepsilon \sqrt{\abs{y_1}}) h_2\paren{\frac{1}{\sqrt{\abs{y_1}}}}, \\
	\wtilde{h_1}(u) :=& \int_0^\infty \e{u\paren{x_2+x_2^{-1}}} h_1\paren{x_2^2} \frac{dx_2}{x_2^2},
\end{align*}
while
\[ \what{F}^d(\mu) = \frac{1}{4} \int_Y T_{w_l}(F)(y) \wbar{K^d_{w_l}(y,\mu)} dy \]
is the long-element Bessel transform of $T_{w_l}(F)$, and $S_{w_l}(m,n,c)$ is the long-element Kloosterman sum.
\end{thm}
Note: The function $T_{w_l}(F)$ is initially given by an integral transform of the test function $F:PGL(3,\R) \to \C$ constructed from $f,h_1,h_2$ in \cref{sect:TestFun} and it reduces to the form above for that particular $F$.
Clearly $T_{w_l}(F)$ and $\what{F}^d(\mu)$ are not uniquely determined by $f$, and we make no claim to have constructed the optimal such test function.

We parameterize the support and derivatives of $f$ and $h_2$ by considering $f$ of the form $f(C^3 \cdot)$ and $h_2$ of the form $h_2(X \cdot)$ where $C,X > 0$.
For $X$ large enough relative to $C$, the long-element term of \cref{thm:MainThm} is zero and we have the spectral expansion of a pure sum of hyper-Kloosterman sums.
Thus \cref{thm:MainCor} is a corollary of \cref{thm:MainThm} by taking $X$ in the range
\begin{align}
\label{eq:XBd}
	X > \frac{C \abs{m_1^2 m_2 n_1 n_2^2}^{1/3}}{200 T_1^{1/3}}.
\end{align}

For $d \ge 2$, we write the spectral parameters of a vector-valued form $\varphi$ of minimal weight $d$ as $\mu_\varphi = \mu^d(r_\varphi) = \paren{\frac{d-1}{2}+r_\varphi,-\frac{d-1}{2}+r_\varphi,-2r_\varphi}$ for some $r_\varphi \in i\R$.
Let $\theta \le \frac{5}{14}$ be the bound towards the Ramanujan-Selberg conjecture for $GL(3)$ \cite[App 2, Prop 1]{KS01}; conjecturally, we should be able to take $\theta=0$.
Then a precise estimate of the spectral side gives
\begin{thm}
\label{thm:MainInDepth}
For $f, h_2$ as in \cref{thm:MainCor} and $C,X \ge 1$,
\begin{align*}
	\sum_{d=0}^\infty \int_{\mathcal{B}_3^{d*}} \frac{\lambda_\Xi(n)\wbar{\lambda_\Xi(m)}}{L(1,\AdSq\Xi)} \what{F}^d(\mu_\Xi) d_H\Xi &= \frac{4}{3} \sum_{d=0}^\infty \sum_{\substack{\varphi \in \mathcal{B}_{3,\text{cusp}}^{d*} \\ \norm{\mu_\varphi} > X^{1/6}}} \frac{\lambda_\varphi(n)\wbar{\lambda_\varphi(m)}}{L(1,\AdSq\varphi)} \wtilde{F}^d_0(\mu_\varphi) \\
	& \qquad + \BigO[f,m,n,\epsilon]{C^{3+\epsilon} X^\epsilon\paren{X^{-\frac{1}{2}}+C^{-\frac{3}{2}}X^{\frac{1}{2}}+C^{3\theta}X^{-1}}},
\end{align*}
where $\wtilde{F}^d_0(\mu)$ defined in \cref{sect:SpectralSide} is supported on the tempered spectrum $\Re(\mu)=0$ or $\Re(r)=0$ and satisfies:
\begin{enumerate}
\item For $d=0,1$, $\wtilde{F}^d_0(\mu)$ is negligible unless $\min_i \abs{\mu_i} \ll X^\epsilon$ in which case we let $w$ be the Weyl group element (i.e. permutation of the coordinates of $\mu$) such that $\mu^w = (-2it_1,it_2-it_1,-it_2-it_1)$ with $t_1 \ll X^\epsilon$ and $t_2 \gg X^{1/6}$ and we have
\[ \wtilde{F}^d_0(\mu) = (-1)^d \frac{i\pi C^3}{2 X t_1} \sgn(w) \paren{\frac{C^3 t_2}{8\pi^3}}^{2it_1} \sum_{\varepsilon \in \set{\pm 1}^2} \chi_d^w((-1,\varepsilon_2)) \varepsilon_1 \what{f}(\varepsilon_2,-1-2it_1) \what{h}_2\paren{\varepsilon_1 \frac{t_2^2}{4\pi^2 X}}. \]

\item For $d \ge 2$, $\wtilde{F}^d_0(\mu^d(r))$ is negligible unless $r \ll X^\epsilon$ in which case for $d \gg X^{1/6}$, we have
\[ \wtilde{F}^d_0(\mu^d(r)) = \frac{\pi i^d C^3 (d+1+2r)}{X d^2} \paren{\frac{C^3}{4\pi^3 d}}^{2r} \sum_{\varepsilon \in \set{\pm 1}^2} \varepsilon_2^d \what{f}(\varepsilon_2,-1-2r) \what{h}_2\paren{\varepsilon_1 \frac{(d-1)^2}{16\pi^2 X}}. \]
\end{enumerate}
Here
\begin{align*}
	\what{f}(\varepsilon,s) :=& \int_0^\infty f(\varepsilon y_2) y_2^s \frac{dy_2}{y_2}, &
	\what{h_2}(\xi) :=& \int_0^\infty \e{-\xi y_1} h_2(y_1) dy_1,
\end{align*}
are the (signed) Mellin transform of $f$ and the one-sided Fourier transform of $h_2$, respectively, and $\chi_d^w$ is a character defined in \cref{sect:Background}.
\end{thm}

The error term $C^{3+3\theta+\epsilon}X^{\epsilon-1}$ is a direct result of some (conjecturally non-existent) low-lying complementary series forms, i.e. having spectral parameters of the shape $\mu_\varphi=(\sigma+it,-\sigma+it,-2it)$ (up to permutation) with $\sigma > 0$ close to $\theta$ and $t \ll X^\epsilon$.
For applications, these terms can, in principle, be handled by writing them out explicitly or by taking $X$ slightly larger than $C^{3\theta}$ (or by assuming a stronger bound toward Ramanujan-Selberg, say $\theta \le \frac{1}{6}$).
Note that density theorems such as \cite[Thm 2]{Val01} do not help us here, since we're looking at the forms with low-lying spectral parameters.

Treating $m,n$ and $f$ as fixed, we see in \cref{sect:LarsenBd,sect:SpectralSide} that the trivial bound (i.e. applying the Deligne bound to $S_{w_5}$) on the hyper-Kloosterman side of \cref{thm:MainCor} is $C^{3+\epsilon}$, while the trivial bound on the spectral side is $C^3 X^{\frac{1}{2}+\epsilon}$.
This behavior is not unexpected:
In one-variable terms, the long element Kloosterman zeta function converges on $\Re(s_i) > \frac{1}{2}$, while the hyper-Kloosterman zeta function converges on $\Re(s_i)>0$, but there are (poles of the Mellin transform of the Whittaker functions of) automorphic forms on $\Re(s_i)=0$, so we win for the long element, but not for the hyper-Kloosterman sums; this can be made precise by phrasing the arithmetic side in terms of the Mellin expansion of the Bessel functions (which have poles on $\Re(s_i)=0$).

We note that the increase in the trivial bound does not rule out practical applications:
\begin{itemize}
\item The purpose of such a formula is to change the nature of the summation in hopes of finding cancellation, and certainly the nature of the spectral sum is different from that of the sum of hyper-Kloosterman sums.

\item At $X \asymp C$, the trivial bound on the spectral side is precisely halfway between the trivial bound and the term-counting bound $C^4$ (i.e. using $\abs{S_{w_5}(m,n,c)} \le c_1 c_2$) on the hyper-Kloosteman side, so applying the above formula is less costly than opening the hyper-Kloosterman sums, provided $m$ and $n$ are not too large.

\item In case the formula is being used inside sums over the coordinates of $m$ and $n$, one might hope to extract asymptotics for the spectral and long-element terms of \cref{thm:MainThm}, whose first terms must then cancel.

\item If one can handle a moderately long sum of the long-element Kloosterman sums, it makes sense to apply \cref{thm:MainThm} with $X$ less than $C$; see \cref{sect:LEAside}.

\item As usual, the analysis of this paper relates to the place at infinity and applies nearly verbatim if one changes the particular discrete group, nebentypus, cusp, etc.
In particular, in the case where we do not have square-root cancellation bounds for the hyper-Kloosterman sums, these results should become highly nontrivial, with the caveat that we are using the arithmetically-weighted Weyl law which also depends (somewhat) on the the strength of such bounds.
\end{itemize}

\section{Further Explorations}
We summarize the remaining problems related to the $GL(3)$ Kuznetsov formulas:
\begin{itemize}
\item For the remaining $GL(3)$ Kloosterman sums, those corresponding to the transpositions $w_2,w_3$ are the classical Kloosterman sums, for which the summation formulas are due to Kuznetsov, while the summation formula for the remaining order-three Weyl element $w_4$ follows immediately from the above by replacing all of the relevant functions with their duals $f^\iota(g) := f(g^\iota)$, $g^\iota := w_l \paren{\trans{g}}^{-1} w_l$ using $w_4=w_5^\iota$.

\item We have only analyzed the case where $m$ and $n$ are small compared to $c$; in particular, we have assumed $C \ge 1$.
It might yield interesting results to study the case $C < 1$, i.e. $m_1 m_2 n_1 > c_1^3$.

\item In \cref{sect:LEAside}, we see that taking $X$ a little smaller than $C$ results in a sum over the moduli of the long-element Kloosterman sums where one moduli is essentially fixed and the other ranges over a region of size $C^2$.
As the long-element Kloosterman sum factors into a (short sum of) product(s) of two $\operatorname{Kl}_2$ Kloosterman sums \cite{Stevens} (see also \cite{KiralNaka}), one might hope that this would yield to $GL(2)$ methods, e.g. the classical Kuznetsov formulas.

\item A consequence of the expansion of \cref{thm:w5BesselExpand} is that one can construct the formula of \cref{thm:MainThm} by adding the spectral Kuznetsov formulas for all $d$ with individually chosen test functions $\what{F}^d(\mu)$, though this is not the approach of the current paper.
Again, the coefficients $\what{F}^d(\mu)$ are not uniquely determined by $f$ and it remains an open question whether one can choose these coefficients such that
\[ \sum_{d \ge 0} \int_{\frak{a}^d_0} \what{F}^d(\mu) \, K^d_w(y,\mu) \sinmu^{d*}(\mu) d\mu \]
is $f(y)$ for $w=w_5$ and zero unless $w=w_5,w_l$ (as we have done) and yet have the trivial bound on both the spectral expansion and long-element Kloosterman term be comparable to the trivial bound on the sum of hyper-Kloosterman sums.
The analysis of \cref{sect:LEAside} suggests a negative answer.
In any case, an explicit determination of the class of functions $\what{F}^d(\mu)$ that make \eqref{eq:MainCor} true for a given $f$ is highly desirable.

\item We note that the support of $\wtilde{F}^d_0(\mu)$ being effectively at $\min_i \abs{\mu_i}=0$ makes it tantalizing to think that a (limit of) properly constructed test function(s) would result in $\wtilde{F}^d_0(\mu)$ supported on the symmetric squares of $SL(2,\Z)$ cusp forms; in that case, the number of such forms with $\norm{\mu} \ll X^{1/2}$ is $\asymp X$ (by the $SL(2,\Z)$ Weyl law) so the trivial bound would become $C^3 X^{-1/2+\epsilon} \asymp C^{5/2+\epsilon}$ at $X\asymp C$.
(This relates to the so-called ``spectral gap'' and ``conductor-dropping'' phenomena in the sense that it suggests there are few, if any, points in the spectrum close to the symmetric squares and the apparent difference in the trivial bounds is due to the relative sparsity of symmetric squares in the whole $SL(3, \Z)$ spectrum.)
Furthermore, this leads one to think that the proper trace formula for hyper-Kloosterman sums should involve an expansion over just these symmetric squares.
Unfortunately, the effect of $f$ on the spectral expansion is precisely orthogonal to the self-dual forms -- i.e. if we write $\mu=(a,-a,0)+(r,r,-2r)$ with $a$ large, then $\wtilde{F}^d_0(\mu)$ will depend on $f$ solely through $\what{f}(\cdot,-1-2r)$ (i.e. there is no interaction between $f$ and $a$, though deeper asymptotics do depend on the moments $\what{f}(\cdot,-1-2r+n)$, $n \in \N$), but the self-dual forms occur at $r=0$.


\item Lastly, one approach to isolating the symmetric squares in the $SL(3,\Z)$ spectrum (\'a la Langlands' Beyond Endoscopy \cite{Langlands03}) would be to take $\what{f}$ approaching a point mass (e.g. taking $f$ of the form $f(y^\delta)$ for $\delta > 0$ small), but this is beyond the scope of the current paper.
\end{itemize}

\section{Background}
\label{sect:Background}
Define $G=PGL(3,\R)$, $K = PO(3,\R)$, $\Gamma = SL(3,\Z)$,
\[ U(R) = \set{\Matrix{1&x_2&x_3\\&1&x_1\\&&1} \setdiv x_i \in R}, \qquad Y = \set{\Matrix{y_1 y_2\\&y_1\\&&1} \setdiv y_i \in \R^\times}, \]
$V = Y \cap K$, $Y=Y^+ V$, with the appropriately normalized Haar measures
\[ dx = dx_1 \, dx_2 \, dx_3, \qquad dy = \frac{dy_1 \, dy_2}{\abs{y_1 y_2}^3}, \qquad \int_K dk = 1, \qquad dg=dx \, dy \, dk. \]
We tend not to distinguish between elements of $Y$ and pairs $(y_1,y_2) \in (\R^\times)^2$ as the multiplication is the same.

On $K$, we apply the coordinates
\[ k=\Matrix{\cos\theta_1&-\sin\theta_1&0\\ \sin\theta_1&\cos\theta_1&0\\0&0&1}\Matrix{\cos\theta_2&0&-\sin\theta_1\\0&1&0\\ \sin\theta_2&0&\cos\theta_1}\Matrix{\cos\theta_3&-\sin\theta_3&0\\ \sin\theta_3&\cos\theta_3&0\\0&0&1}, \]
with $\theta_1,\theta_3 \in [0,2\pi)$, $\theta_2 \in [0,\pi]$, and Haar probability measure is
\[ dk = \frac{1}{8\pi^2} \sin\theta_2 \, d\theta_1 \, d\theta_2 \, d\theta_3. \]

Define characters of $U(\R)$ by $\psi_y(x) = \psi_I(yxy^{-1}) = \e{y_1 x_1+y_2 x_2}$, $\e{x} = e^{2\pi i x}$ and of $Y$ by
\[ p_\mu\Matrix{a_1\\&a_2\\&&a_3} = \prod_{i=1}^3 \abs{a_i}^{\mu_i}, \qquad \chi_\delta\Matrix{a_1\\&a_2\\&&a_3} = \prod_{i=1}^3 \sgn(a_i)^{\delta_i}, \]
where $\mu \in \C^3, \delta \in \Z^3$ satisfy $\mu_1+\mu_2+\mu_3=0$, $\delta_1+\delta_2+\delta_3 \equiv 0 \pmod{2}$.
Call $\mu$ and $\delta$ the (Harish-Chandra/Langlands) spectral parameters of an automorphic representation/automorphic form.

The representations of $K$ are given by the Wigner-$\mathcal{D}$ matrices.
Up to equivalence, there is one (irreducible, unitary) representation $\mathcal{D}^d:K\to GL(2d+1,\C)$ of each odd dimension $2d+1$, $d \ge 0$, where $d=0$ corresponds to the trivial representation.

The Weyl group in $G$ is, up to scalar multiples, the set of all matrices with exactly one 1 on each row and column.
We write $\Weyl=\set{I,w_2,w_3,w_4,w_5,w_l} = \Weyl_3 \cup w_2\Weyl_3$ with $\Weyl_3=\set{I, w_4, w_5}$, where
\begin{gather*}
	w_2=-\Matrix{&1\\1\\&&1}, \qquad w_3=-\Matrix{1\\&&1\\&1}, \\
	w_4=\Matrix{&1\\&&1\\1}, \qquad w_5=\Matrix{&&1\\1\\&1}, \qquad w_l=-\Matrix{&&1\\&1\\1}.
\end{gather*}
The Weyl group acts on $\mu$ and $\delta$ by $p_{\mu^w}(y)=p_\mu(w y w^{-1})$ and $\chi_{\delta^w}(y) = \chi_\delta(wyw^{-1})$.
We also define $\wbar{U}_w=(w^{-1} \trans{U} w) \cap U$, $U_w=(w^{-1} U w) \cap U$ and
\[ Y_w = \set{y \in Y \setdiv \psi_y(wuw^{-1})=\psi_I(u), \forall u \in U_w(\R)}. \]
In particular, $Y_{w_5} = \set{(1,y_2) \in Y}$.

The Jacquet-Whittaker function is given by the (Jacquet) integral
\[ W^d(g,\mu,\delta,\psi) = \int_{U(\R)} I^d_{\mu,\delta}(w_l x g) \wbar{\psi(x)} dx, \]
where
\[ I^d_{\mu,\delta}(xyk) = p_{\rho+\mu}(y) \Sigma^d_\delta \mathcal{D}^d(k), \qquad \Sigma^d_\delta = \tfrac{1}{4} \sum_{v \in V} \chi_\delta(v) \mathcal{D}^d(v). \]
Since $I^d_{\mu,\delta}$ is an eigenfunction of the Casimir operators, so is the Whittaker function, with the same eigenvalues.
For convenience, we write $W^d(g,\mu) := W^d(g,\mu,(d,d,0),\psi_I)$.

\subsection{The Weyl law}
\label{sect:WeylLaw}
For $d=0,1$, define $\frak{a}^d_0=\set{\mu \setdiv \Re(\mu)=0}$ and for $d \ge 2$, define $\frak{a}^d_0 = \set{\mu^d(r) \setdiv \Re(r)=0}$.
We refer to such $\mu$ as ``tempered''.

We say $\mu$ of the form $\mu=(x+it,-x+it,-2it)$ (up to permutation) with $0<x<\frac{1}{2}$ is ``complementary''.
It is known that the spectral parameters of automorphic forms for $SL(3,\Z)$ are either tempered or complementary with $x < \theta$ where $\theta \le \frac{5}{14}$ is the Kim-Sarnak bound toward the Ramanujan-Selberg conjecture.

The arithmetically-weighted Weyl laws for $SL(3,\Z)$ cusp forms are given in \cite{Val01,WeylI,WeylII}, but we need a weaker bound for some oddly-shaped regions, so we note that taking a test function of the form
\[ H(\mu) = -\int_\Omega \exp \sum_{i=1}^3 (\mu_i-\mu'_i)^2 d\mu' \]
in \cite[Thm 9]{WeylI} easily implies for $d=1$ that
\begin{align}
\label{eq:WeylLaw01}
	\sum_{\substack{\varphi \in \mathcal{B}^{d*}_{3,\text{cusp}}\\ \mu_\varphi\in\Omega}} \frac{1}{L(1,\AdSq \varphi)} \ll \int_{\Omega+B(T^\epsilon,0)} \prod_{i<j} \paren{1+\abs{\mu_i-\mu_j}} \abs{d\mu},
\end{align}
where $\Omega \subset \frak{a}^d_0$, $T \ge 1$ is the diameter of $\Omega$ and $B(r,0) \subset \frak{a}^d_0$ is the ball of radius $r$ centered at $\mu=0$.

Note that the Bessel functions for $d=0$ differ from those of $d=1$ only in the particular linear combinations of the Mellin-Barnes integrals, but the proof of \cite[Thm 9]{WeylI} applied the triangle inequality to that sum, so \eqref{eq:WeylLaw01} applies equally well to $d=0$.
One might worry, after reading \cref{sect:FormsNearWalls}, that the paper \cite{WeylI} failed to handle the case with $\abs{\mu_1-\mu_2}$ small (and this certainly was an oversight by the author), but in that situation this corresponds to a zero of the spectral measure (which is not present in the current situation); the same is true for the $d=0$ case.

For $d \ge 2$ and $T \ge 1$, \cite[Thm 1]{WeylII} implies
\begin{align}
\label{eq:WeylLaw2plus}
	\sum_{\substack{\varphi \in \mathcal{B}^{d*}_{3,\text{cusp}}\\ \abs{r_\varphi} \le T}} \frac{1}{L(1,\AdSq \varphi)} \ll dT(d+T)^2.
\end{align}
Note:
There is a slight caveat that \cite[Thm 1]{WeylII} only applies to $d \ge 3$, but instead \cite[Prop 12]{WeylII} implies the weaker bound $T^4$ for $d=2$.
The final paragraph of \cite[Sect 3.3]{ArithKuzII} describes how the Weyl law for $d=2$ can be obtained; alternately, we can simply use the weaker bound here, as the large contribution to \cref{thm:MainInDepth} on $d \ge 2$ comes from $r$ small with $d$ large.

Lastly, the complementary spectrum for $d=0$ can be handled by a comment at the bottom of page 678 in \cite{Val01} and for $d=1$ by a comment in the first paragraph on page 27 of \cite{WeylI}; in both cases, for $T \ge 1$, we have
\begin{align}
\label{eq:WeylLawCompSer}
	\sum_{\substack{\varphi \in \mathcal{B}^{d*}_{3,\text{cusp}}\\ \mu_\varphi \text{ complementary} \\ \abs{t_\varphi} \le T}} \frac{1}{L(1,\AdSq \varphi)} \ll T^{3+\epsilon}.
\end{align}

\subsection{Mellin-Barnes integrals of the $GL(3)$ Bessel functions}
\label{sect:BesselMBIntegrals}
Define the (normalized) Mellin transform of the Bessel functions by
\begin{equation}
\label{eq:KwlMellinTrans}
\begin{aligned}
	K^d_{w_l}(y,\mu) =:& \int_{-i\infty}^{+i\infty} \int_{-i\infty}^{+i\infty} \abs{4\pi^2 y_1}^{1-s_1} \abs{4\pi^2 y_2}^{1-s_2} \what{K}^d_{w_l}(s,\sgn(y),\mu) \frac{ds_1}{2\pi i} \frac{ds_2}{2\pi i}, \\
	\what{K}^d_{w_l}(s,v,\mu) =& \frac{1}{(2\pi)^8} \int_{Y^+} K^d_{w_l}(vy,\mu) (4\pi^2 y_1)^{s_1+1} (4\pi^2 y_2)^{s_2+1} dy.
\end{aligned}
\end{equation}
The unbounded portion of the contours in each $\int_{-i\infty}^{i\infty} \ldots ds_i$ must pass to the left of $\Re(s_i)=0$ (and the finite part must pass to the right of the poles of the integrands) to maintain absolute convergence for $\sgn(y) \ne (+,+)$.

Write $v = \sgn(y)$, $\chi_d^w = \chi_{(d,d,0)^w}$ and
\[ \frac{1}{\cscmu^d(\mu)} := 8 \sin\tfrac{\pi}{2}(\mu_1-\mu_2) \sin\tfrac{\pi}{2}(d+\mu_1-\mu_3) \sin\tfrac{\pi}{2}(d+\mu_2-\mu_3), \]
then for $d=0,1$, we have
\begin{equation}
\label{eq:Kwl1Sgn}
\begin{aligned}
	\what{K}^d_{w_l}(s,v,\mu) =& \frac{1}{2 \pi^2} \cscmu^d(\mu) \sum_{w \in \Weyl_3} \chi_d^w(v) \, G^{v_1,v_2}(s,\mu^w).
\end{aligned}
\end{equation}
where
\begin{align*}
	G^{++}(s,\mu) :=& \frac{\Gamma(s_1-\mu_1) \Gamma(s_1-\mu_2) \Gamma(s_1-\mu_3) \Gamma(s_2+\mu_1) \Gamma(s_2+\mu_2) \Gamma(s_2+\mu_3)}{3\pi^2 \Gamma(s_1+s_2)} \prod_{i<j} \sin \pi(\mu_i-\mu_j), \\
	G^{+-}(s,\mu) :=& \sin\pi(\mu_2-\mu_3) \frac{\Gamma(s_1-\mu_2) \Gamma(s_1-\mu_3) \Gamma(s_2+\mu_1) \Gamma(1-s_1-s_2)}{\Gamma(1-s_1+\mu_1) \Gamma(1-s_2-\mu_2) \Gamma(1-s_2-\mu_3)}, \\
	G^{-+}(s,\mu) :=& \sin\pi(\mu_1-\mu_2) \frac{\Gamma(s_1-\mu_3) \Gamma(s_2+\mu_1) \Gamma(s_2+\mu_2) \Gamma(1-s_1-s_2)}{\Gamma(1-s_1+\mu_1) \Gamma(1-s_1+\mu_2) \Gamma(1-s_2-\mu_3)}, \\
	G^{--}(s,\mu) :=& \sin\pi(\mu_1-\mu_3) \frac{\Gamma(s_1-\mu_1) \Gamma(s_1-\mu_3) \Gamma(s_2+\mu_1) \Gamma(s_2+\mu_3)}{\Gamma(1-s_1+\mu_2) \Gamma(1-s_2-\mu_2) \Gamma(s_1+s_2)}.
\end{align*}
Note: Compared to \cite[Sect 6.14.3]{ArithKuzII}, we have factored out a term $1/\prod_{i<j} \sin \pi(\mu_i-\mu_j)$ from the $\Weyl_3$ sum, which is invariant under $\Weyl_3$ (since it consists of even permutations).

For $d \ge 2$, it is somewhat more useful to write these as
\begin{align}
\label{eq:KwldMB}
	\what{K}_{w_l}^d(s,v,\mu^d(r)) =& \frac{1}{4\pi^2} (-\varepsilon_1\varepsilon_2)^d B^{\varepsilon_1,\varepsilon_2}_{w_l}\paren{s,r} Q(d,s_1-r) Q(d,s_2+r),
\end{align}
where
\begin{align*}
	Q(d,s) :=& \frac{\Gamma\paren{\frac{d-1}{2}+s}}{\Gamma\paren{\frac{d+1}{2}-s}}, &
	B^{\varepsilon_1,\varepsilon_2}_{w_l}(s,r) := \piecewise{
	0& \If \varepsilon=(+,+), \\
	B\paren{s_1+2r,1-s_1-s_2}& \If \varepsilon=(+,-), \\
	B\paren{s_2-2r,1-s_1-s_2}& \If \varepsilon=(-,+), \\
	B\paren{s_1+2r,s_2-2r}& \If \varepsilon=(-,-),}
\end{align*}
and $B(a,b)$ is the usual beta function.

\subsection{Oscillatory integrals}
We will engage in some stationary phase analysis of oscillatory integrals.
Suppose $\omega:\R\to\C$ is smooth with fixed compact support and satisfies $\omega^{(j)}(t) \ll_j T U^{-j}$ for some $T,U > 0$ and all $j \ge 0$.
Also suppose $\phi:\R\to\R$ is smooth and satisfies $\phi^{(j)}(t) \ll_j Y Q^{-j}$ on the support of $\omega$ for some $Y, Q > 0$ and all $j \ge 2$.
Then we consider the integral
\[ \mathcal{I} := \int_{-\infty}^\infty \omega(t) \e{\phi(t)} dt. \]

\begin{lem}[{\cite[Lem 8.1]{BKY01}}]
\label{lem:BKY8.1}
Suppose that $\abs{\phi'(t)} \ge R$ on the support of $\omega$ for some $R > 0$, then for any $A>0$,
\[ \mathcal{I} \ll_A T\paren{\paren{\frac{QR}{\sqrt{Y}}}^{-A}+(RU)^{-A}}. \]
In particular, for any $B > 1$, $\mathcal{I} \ll_{\epsilon,A} T B^{-A}$, unless $R < B^\epsilon \paren{\sqrt{Y}/Q+1/U}$.
\end{lem}

\begin{prop}[{\cite[Prop 8.2]{BKY01}}]
\label{prop:BKY8.2}
Suppose that there is a unique $t_0$ in the support of $\omega$ such that $\phi'(t_0)=0$, that $\omega$ is supported on an interval of length $\ll U_1 \ge U$, that $\phi''(t) \gg Y Q^{-2}$ on the support of $\omega$ and that for some $\delta > 0$, $Y \ge Z^{3\delta}$ and $U \ge \frac{Q Z^{\delta/2}}{Y^{1/2}}$ with $Z := Q+T+Y+U_1+1$, then for any $A > 0$,
\[ \mathcal{I} = \frac{\e{\phi(t_0)}}{\sqrt{\abs{\phi''(t_0)}}} \sum_{0 \le j \le 3A/\delta} p_j(t_0)+\BigO[A,\delta]{Z^{-A}}, \]
where
\[ p_j(t) := \frac{\e{\frac{1}{8}\sgn(\phi''(t))}}{j!} \paren{\frac{i}{4\pi\phi''(t_0)} \frac{d^2}{dt^2}}^j \omega(t) \e{\phi(t)-\phi(t_0)-\tfrac{1}{2}\phi''(t_0)\paren{t-t_0}^2} \]
satisfies
\[ \frac{d^k}{dt_0^k} p_j(t_0) \ll_{j,k} T \paren{U^{-k} + Q^{-k}}\paren{\paren{U^2 Y/Q^2}^{-j}+Y^{-j/3}}. \]
\end{prop}

\section{Preliminaries on Bessel functions}
In \cref{sect:BesselBd}, we will show
\begin{lem}
\label{lem:BesselBd}
For $\sigma \ge 0$, $t \in \R$, $x >0$,
\begin{align}
\label{eq:PhragLindJBd}
	J_{\sigma+it}(x) \ll (1+\abs{x^2-\sigma^2+t^2}+\sigma \abs{t})^{-1/4} e^{\frac{\pi}{2}\abs{t}}.
\end{align}
In particular, for $-1 \le \sigma < x/2$ and $\abs{t} < x$, we have
\begin{align}
\label{eq:JsigmaitBd}
	J_{\sigma+it}(x) \ll x^{-1/2} e^{\frac{\pi}{2}\abs{t}}.
\end{align}
Similarly, for $-1 \le \sigma < x/2$ and $x \ge \Max{\abs{t},1}$, we have
\begin{align}
\label{eq:YsigmaitBd}
	Y_{\sigma+it}(x) \ll x^{-1/2} e^{\frac{\pi}{2}\abs{t}}.
\end{align}
\end{lem}
These bounds are simple consequences of the known asymptotic expansions.

In \cref{sect:BesselDervBd}, we will show
\begin{lem}
\label{lem:BesselDervBd}
For $\sigma \ge 0$, $t \in \R$, $x > 0$, we have
\[ \left.\frac{\partial}{\partial \nu} J_\nu(x) \right|_{\nu=\sigma+it} \ll_\epsilon \paren{x+x^{-1}+\abs{t}}^\epsilon (1+x^2+t^2)^{-1/4} e^{\frac{\pi}{2}t}. \]
\end{lem}
This bound is actually slightly stronger than we need, but we feel a little extra work is warranted as it approaches the expected size (which would replace the epsilon power with a logarithm) of the order derivative.

Define
\begin{align*}
	G_1^-(s_1,\mu) =& -\frac{\Gamma(s_1-\mu_3)}{2\pi \Gamma(1+\mu_1-\mu_2) \Gamma(1+\mu_1-\mu_3) \Gamma(1-s_1+\mu_2)}.
\end{align*}
The following two propositions follow from Stirling's formula, and we leave the proof to the reader:
\begin{prop}
\label{prop:hatKBd01}
Consider $d=0,1$ with $\Re(s)$ and $\Re(\mu)$ in a fixed compact set such that the arguments of the gamma functions in $\what{K}^d_{w_l}(s,v,\mu)$ are bounded away from $-\N_0$.
Suppose also that $\abs{\mu_i - \mu_j} \gg 1$ for all $i \ne j$.
Then we have
\[ \what{K}^d_{w_l}(s,v,\mu) \ll \abs{s_1+s_2}^{\frac{1}{2}-\Re(s_1+s_2)} \prod_{i=1}^3 \abs{s_1-\mu_i}^{\Re(s_1-\mu_i)-\frac{1}{2}} \abs{s_2+\mu_i}^{\Re(s_2+\mu_i)-\frac{1}{2}}, \]
and for $w \in \Weyl$ with $v_1=-1$, we have
\begin{align*}
	\res_{s_2=-\mu^w_1} \what{K}^d_{w_l}((-1,v_2), s,\mu) =& \sgn(w) \cscmu^d(\mu) \paren{\chi_d^w(v) G_1^-(s_1,\mu^w)-\chi_d^{ww_3}(v) G_1^-(s_1,\mu^{ww_3})} \\
	\ll& \abs{\mu^w_2-\mu^w_1}^{\Re(\mu^w_2-\mu^w_1)-\frac{1}{2}} \abs{\mu^w_3-\mu^w_1}^{\Re(\mu^w_3-\mu^w_1)-\frac{1}{2}} \\
	& \qquad \times \abs{s_1-\mu^w_2}^{\Re(s_1-\mu^w_2)-\frac{1}{2}} \abs{s_1-\mu^w_3}^{\Re(s_1-\mu^w_3)-\frac{1}{2}}.
\end{align*}
\end{prop}

\begin{prop}
\label{prop:hatKBd2plus}
Consider $d\ge 2$ with $\Re(s)$ and $\Re(r)$ in a fixed compact set such that the arguments of the gamma functions in $\what{K}^d_{w_l}(s,v,\mu^d(r))$ are bounded away from $-\N_0$.
Then we have
\begin{align*}
	\what{K}^d_{w_l}(s,v,\mu^d(r)) \ll& \abs{s_1+s_2}^{\frac{1}{2}-\Re(s_1+s_2)} \abs{s_1+2r}^{\Re(s_1+2r)-\frac{1}{2}} \abs{s_2-2r}^{\Re(s_2-2r)-\frac{1}{2}} \\
	& \qquad \times \paren{d+\abs{s_1-r}}^{2\Re(s_1-r)-1} \paren{d+\abs{s_2+r}}^{2\Re(s_2+r)-1},
\end{align*}
and for $w \in \Weyl$ with $v_1=-1$, we have
\begin{align*}
	\res_{s_2=2r} \what{K}^d_{w_l}(v, s,\mu^d(r)) =& \frac{v_2^d}{4\pi^2} Q(d,-r) Q(d,s_1-r) \\
	\ll& \paren{d+\abs{r}}^{2\Re(r)-1} \paren{d+\abs{s_1-r}}^{2\Re(s_1-r)-1}.
\end{align*}
\end{prop}
Note that, in particular, the exponential parts of Stirling's formula at worst cancel, see \cite[Sect 3.3]{SubConv}.

\section{$w_5$ Bessel functions and Kloosterman sums}
\subsection{Bessel expansions}
For $i=1,2$, let $\Delta_i$ be the degree-$(i+1)$ Casimir operator on $G$.
\begin{prop}
\label{prop:Density}
Let $F(g)$ satisfying $F(ug) = \psi_I(u) F(g)$ be smooth and compactly supported on $U(\R)\backslash G \cong YK$.
For $w \in W$ and $y \in Y_w$, define
\[ T_w(F)(y) := \int_{\wbar{U}_w(\R)} F(ywu) \wbar{\psi_I(u)} du, \]
then for any $w \in W$, we have 
\[ T_w(F)(y) = \sum_{d \ge 0} \int_{\frak{a}^d_0} \what{F}^d(\mu) \, K^d_w(y,\mu) \sinmu^{d*}(\mu) d\mu \]
where
\[ \what{F}^d(\mu) = \frac{1}{4} \int_Y T_{w_l}(F)(y) \wbar{K^d_{w_l}(y,\mu)} dy, \]
and we have
\[ T_{w_l}(F)(y) \ll \abs{y_1 y_2}^{1+\frac{1}{8}}, \qquad (-1)^{(i+1)j} \lambda_i(\wbar{\mu})^j \what{F}^d(\mu) = \what{\Delta_i^j F}^d(\mu).  \]
\end{prop}
Note: The $(-1)^{(i+1)j}$ is because $\Delta_2$ is antisymmetric.

Before we proceed with the proof, some discussion of the Wallach--van den Ban Whittaker inversion formula is required.
The Whittaker-Plancherel formula (aka Whittaker inversion) is an isomorphism of two Hilbert spaces (which are space-consuming to describe) and was originally due to Wallach \cite{Wallach}, but some gaps in the proof were discovered, particularly in \cite{vandenBanKuit}.
This gap is closed by \cite{WallachCorrect} in many cases (in particular for $PSL(3,\R)$) and an alternate argument by van den Ban should be forthcoming.
It is worth mentioning that the case for $PSL(3,\R)$ is very nearly proved in \cite{ArithKuzII} and its preceeding papers:
There are two parts to the theorem -- the Whittaker expansion of a function and the orthogonality of the Whittaker functions (compare the proof below).
The orthogonality at the minimal $K$-types is proved in \cite{GoldKont, WeylI, WeylII} and the action of the Lie algebra implies the orthogonality holds away from the minimal $K$-types.
The Whittaker expansion follows from Harish-Chandra's Plancherel theorem and the fact that the Fourier transform of a zonal spherical function is a product of two Whittaker functions, see \cite[Prop 29]{ArithKuzII} (the use of the minimal-weight vector there is unnecessary).
Unfortunately, this last fact is only proved for $d=0,1$ there; see the discussion in \cite[Sect 9.6]{ArithKuzII}.

\begin{proof}
The Whittaker inversion theorem as in \cite[Thm 10]{ArithKuzII} (i.e. the Whittaker expansion) tells us that
\begin{align}
\label{eq:WhittExpand}
	F(g) =& \sum_{d \ge 0} \int_{\frak{a}^d_0} F^d(g,\mu) \sinmu^{d*}(\mu) d\mu,
\end{align}
where
\[ \what{F}^d(g',\mu) = \sum_{d' \ge 0} (2d'+1) \Tr\paren{W^{d'}(g',\mu) \int_{U(\R)\backslash G} F(g) \trans{\wbar{W^{d'}(g,\mu)}} dg} \]
has super-polynomial decay in $\norm{\mu}$ (by \cite[Lem 9]{ArithKuzII} and the bounds of \cite[Cor 21]{ArithKuzII}) and satisfies
\[ (-1)^{(i+1)j} \lambda_i(\wbar{\mu})^j \what{F}^d(g,\mu) = \what{\Delta_i^j F}^d(g,\mu). \]

For a sufficiently nice (e.g. holomorphic on a tube containing $\frak{a}^d_0$ and rapid decay in $\Im(\mu)$) test function $H(\mu)$, the definition of the Bessel functions (see \cite[(79)]{ArithKuzII}) is
\[ T_w\paren{\int_{\frak{a}^d_0} H(\mu) W^{d'}(\cdot g', \mu) d\mu}(g) = \int_{\frak{a}^d_0} H(\mu) K^d_w(g,\mu) W^{d'}(g', \mu) d\mu. \]
Again, \cite[Lem 9 and Cor 21]{ArithKuzII} show that
\[ H(\mu) = \sinmu^{d*}(\mu) \int_{U(\R)\backslash G} F(g) W^{d'}(g,\mu) dg \]
is sufficiently nice (we have used the fact that $-\wbar{\mu} \in \set{\mu,\mu^{w_2}}$ to replace $\trans{\wbar{W^{d'}(g,\mu)}}$ with $W^{d'}(g,\mu)$) and also justify reordering the $T_w$ integral and the $d$ and $d'$ sums, which tells us that
\begin{align}
\label{eq:TwExpand1}
	T_w(F)(y) = \sum_{d \ge 0} \int_{\frak{a}^d_0} \what{F}^d(I,\mu) K^d_w(y,\mu) \sinmu^{d*}(\mu) d_\mu.
\end{align}
(Note that $\what{F}^d(I,\mu)$ is independent of $w$.)

In particular, we have
\[ T_{w_l}(F)(y) = \sum_{d \ge 0} \int_{\frak{a}^d_0} \what{F}^d(I,\mu) K^d_{w_l}(y,\mu) \sinmu^{d*}(\mu) d_\mu, \]
but in \cite[Lem 27]{ArithKuzII}, we proved orthogonality of the Bessel functions to get
\[ \frac{1}{4} \int_Y T_{w_l}(F)(y) \wbar{K^d_{w_l}(y,\mu)} dy = \what{F}^d(I,-\wbar{\mu}). \]
As in the comment preceding that lemma, $\what{F}^d(I,-\wbar{\mu}) = \what{F}^d(I,\mu)$ whenever $-\wbar{\mu} \in \set{\mu,\mu^{w_2}}$.
The proposition then follows from \eqref{eq:TwExpand1}, noting the bound in \cite[Lem 25]{ArithKuzII}.
(Again, the interchange of the $y$-integral with the $d$ sum follows by rapid decay of $\what{F}^d(I,\mu)$.)
\end{proof}
For the test function $F$ constructed in \cref{sect:TestFun}, it will turn out that $T_{w_5}(F)(y) = f(y)$ while $T_w(F)(y) = 0$ unless $w=w_5,w_l$.

\subsubsection{Aside on some worrisome but ultimately irrelevant questions}
\begin{itemize}
\item Note that the iterated integral $\what{F}^d(\mu)$ can be written in the form
\[ \int_{U(\R)\backslash G} F(g) \wbar{K^d_{w_l}(gw_l,\mu)} dg. \]
This is slightly non-trivial:
In the conversion $utw_lu'=xyk$ between Bruhat and Iwasawa coordinates on the long-element Weyl cell, the function $\frac{K^d_{w_l}(t,\mu)}{\abs{t_1 t_2}}$ is (relatively) well-behaved in both coordinates by the bound \cite[Lem 22]{ArithKuzII}, but the factor $\abs{t_1 t_2}=y_1 y_2 \abs{\sec \theta_1 \csc^2 \theta_2 \sec \theta_3}$ has singularities on $K$, though the factor $\sin\theta_2$ in the measure $dk$ offsets one factor $\csc \theta_2$.
(The measures compare as $\frac{1}{4} du \, dt \, du' = dx \, dy \, dk$.)
So some care has to be taken with integrals against the $K^d_{w_l}$ function, and in this integral representation of $\what{F}^d(\mu)$, we are actually using the slight decay of the factor $\frac{K^d_{w_l}(t,\mu)}{\abs{t_1 t_2}}$ at $\cos\theta_1 \sin\theta_2 \cos\theta_3=0$ (i.e. $\abs{t_1 t_2}=\infty$) to get convergence.

\item We have
\[ T_{w_l}(\Delta_i^j F) = \wtilde{\Delta}_i^j T_{w_l}(F) \]
using the operators $\wtilde{\Delta}_i$ in \cite[Prop 3]{ArithKuzII}, but it appears that $\wtilde{\Delta}_i^j T_{w_l}(F)$ has growth like $y_1^{c+(i+1)j}$ while $T_{w_l}(\Delta_i^j F)$ does not; see \cref{sect:BruhatLikeLie}.
The apparent discrepancy can be resolved by writing out
\[ 0 = \int_{U(\R)} \frac{\partial}{\partial x_1} F(yw_l x) \wbar{\psi_I(x)} dx = \int_{U(\R)} \frac{\partial}{\partial x_2} F(yw_l x) \wbar{\psi_I(x)} dx \]
and comparing these expressions to the terms of $0=T_{w_l}(\Delta_1 F) - \wtilde{\Delta}_1 T_{w_l}(F)$, etc.
One can also see that this does not imply rapid decay in $y_1$ since multiple terms involving a factor $y_1$ in $\wtilde{\Delta}_1 T_{w_l}(F)$ are required to trade for the corresponding terms in $T_{w_l}(\Delta_1 F)$.

\item So we have that
\[ \int_Y T_{w_l}(F)(y) \wbar{K^d_{w_l}(y,\mu)} dy \]
decays rapidly in $\norm{\mu}$ and $d$, and we are left to wonder why we cannot apply the arithmetic Kuznetsov formula to the long-element term in the hyper-Kloosterman Kuznetsov formula to conclude that the sum of hyper-Kloosterman sums is simply zero.
The answer is that taking $\wtilde{f}(y) = T_{w_l}(F)(y)$ in the arithmetic Kuznetsov formula means the corresponding $\wtilde{F}(g)$ is not bounded much less smooth or even continuous on $K$, and hence the spectral expansion fails.
(Particularly, the expansion over the spectral basis of $L^2(K)$ converges too slowly.)
One can correct this by inserting a smooth partition of unity to enforce $\abs{y_1} \le X$ in $\wtilde{f}(y)$, then applying the arithmetic Kuznetsov formula, and finally taking the limit as $X \to \infty$.
This is effectively what we are doing except we don't involve the arithmetic Kuznetsov formula (no need) and we don't take the limit as $X \to \infty$ (doesn't converge in any obvious manner, aside from the $X$-invariance of the left-hand side in \cref{thm:MainCor}).
\end{itemize}

\subsection{Poincar\'e series}
\label{sect:FormalPoin}
\begin{prop}
\label{prop:SpectralInterpret}
For $m,n\in \Z^2$ with $m_1m_2n_1n_2\ne0$ and $F(g)$ satisfying $F(ug) = \psi_I(u) F(g)$ and smooth and compactly supported on $U(\R)\backslash G$, we have
\begin{align*}
	\sum_{w \in W} \sum_{c \in \Z_{\ne0}^2} S_w(m,n,c) p_{\rho^w}(n) T_w(F)(mcwn^{-1}w^{-1}) &= p_\rho(m) \sum_{d=0}^\infty \int_{\mathcal{B}_3^{d*}} \frac{\lambda_\Xi(n)\wbar{\lambda_\Xi(m)}}{L(1,\AdSq\Xi)} \what{F}^d(\mu_\Xi) d_H\Xi,
\end{align*}
with $T_w(F)$ and $\what{F}^d(\mu)$ as in \cref{prop:Density}.
\end{prop}
Note: See \cite[Thm 8, (52)]{ArithKuzII} for the details of the terms of the spectral expansion; the $d$ sum and the $\Xi$ integral there are running through the minimal-weight forms, while $d'$ is running through all lifts to higher weights.
In other words, the combined $d$ and $\Xi$ sum/integral is the decomposition of the right regular representation into irreducible representations (organized by their minimal weights), while the $d'$ sum represents the expansion into a basis (via $K$-types) of a given irreducible representation.
We are also including in the measure $d_H\Xi$ the constants associated with the conversion to Hecke eigenvalues; see \cite[Sect 6.7]{ArithKuzII}.

\begin{proof}
The proof follows that of \cite[Sect 8]{ArithKuzII}.

Construct a Poincar\'e series
\[ P_m(g) = \sum_{\gamma\in U(\Z)\backslash\Gamma} F(m\gamma g), \]
(treating $m$ as an element of $Y$) and consider its Fourier coefficient
\[ \mathcal{P} = \int_{U(\Z) \backslash U(\R)} P_m(u n^{-1}) \wbar{\psi_n(u)} du. \]

Then the Bruhat decomposition (see \cite[Sects 2.2.3, 2.2.5]{MeThesis}; the factor $p_{\rho^w-\rho}(n)$ here is the result of the substitution $un^{-1} \mapsto n^{-1} u$ on $\wbar{U}_w(\R)$) gives
\[ \mathcal{P} = \sum_{w \in W} \sum_{c \in \Z_{\ne0}^2} p_{\rho^w-\rho}(n) S_w(m,n,c) T_w(F)(mcwn^{-1}w^{-1}), \]
and the spectral decomposition \cite[Thm 8 as in eq. (52)]{ArithKuzII} and the proof of \cref{prop:Density} (compare \cite[(114)-(116), Lem 27]{ArithKuzII}) give
\[ \mathcal{P} = p_\rho(mn^{-1}) \sum_{d=0}^\infty \int_{\mathcal{B}_3^{d*}} \frac{\lambda_\Xi(n)\wbar{\lambda_\Xi(m)}}{L(1,\AdSq\Xi)} \what{F}^d(\mu_\Xi) d_H\Xi, \]
using the conversion to Hecke eigenvalues as in \cite[Sect 6.7]{ArithKuzII}.
\end{proof}

\subsubsection{Aside on convergence}
In the proof of \cite[Thm 1.1]{HWI} in section 4.5 of that paper, we see that the expansion of the minimal parabolic Eisenstein term involves Mellin inversion in the region of absolute convergence of the Eisenstein series.
Thus, for the expansion of the $1,1,1$ constant term, call it $P^{111}_m(g)$, of $P_m(g)$, it is sufficient to have absolute convergence of
\begin{align*}
	&\int_{Y^+} \norm{\int_K P_m^{111}(yk) \trans{\wbar{\mathcal{D}^d(k)}} dk} p_{(2+\epsilon)\rho}(y) dy \\
	& \ll_\epsilon \int_{Y^+} \int_K \abs{F(myk)} dk \, p_{(2+\epsilon)\rho}(y) \prod_{i=1}^2 (1+y_i^{-1})^3 dy < \infty
\end{align*}
for $\epsilon$ in a subinterval of $(0,\infty)$, using \cite[(4.2)]{HWI}.

Convergence of the cuspidal expansion requires
\begin{enumerate}
\item Square-integrability of $P_m(g)$; this follows as above from $\abs{F(mxyk)} \ll p_{(2+\epsilon)\rho}(y)$ (by compact support) and unfolding one of the $P_m(g)$.

\item Square-integrability of the Eisenstein terms; see the last sentence of \cite[Sect 4.5]{HWI}.
\end{enumerate}
Convergence of the expansions of the $2,1$ and $1,2$ constant terms is intermediate between the cuspidal and $1,1,1$ terms (i.e. requiring Mellin inversion in one variable and then square-integrability in the appropriate copy of $SL(2,\Z)\backslash PSL(2,\R)$), and the proof is similar.

Pointwise convergence of the spectral expansion can be shown via
\[ \abs{\int_{\Gamma\backslash G} P_m(g) \wbar{\Xi^{d'}(g)} dg} \le \int_{U(\Z)\backslash G} \abs{F(mg) \Xi^{d'}(g)} dg \]
from the symmetry of the Casimir operators (for rapid decay in $d'$ and $\mu_\Xi$ and hence also $d$) and a crude bound on the supremum norm of the (appropriately normalized) $GL(3)$ spectral basis elements, say $\Xi^{d'}(xyk) \ll (\norm{\mu_\Xi} + d' + y_1 + y_1^{-1} + y_2 + y_2^{-1})^{100}$.
The latter follows easily from the Fourier expansion \cite[(3.30)-(3.31)]{HWI} and the bounds of \cite[Sect 7.1]{ArithKuzII} and \cite[Thm 2]{XLi01}; see the proof of \cite[Thm 1]{BHM01}.

\subsection{Bruhat-like coordinates}
For a matrix
\[ g = \Matrix{g_{11}&g_{12}&g_{13}\\g_{21}&g_{22}&g_{23}\\g_{31}&g_{32}&g_{33}} \in GL(3,\R), \]
we define the Pl\"ucker coordinates of $U(\R)\backslash GL(3,\R)$ as
\begin{align*}
	A_1(g)=&g_{31}, & B_1(g)=&g_{32}, & C_1(g) =& g_{33}, \\
	A_2(g)=&g_{21}g_{32}-g_{22}g_{31}, & B_2(g)=&g_{23}g_{31}-g_{21}g_{33}, & C_2(g) =& g_{22} g_{33}-g_{23}g_{32},
\end{align*}
and when $B_1(g),A_2(g) \ne 0$, we may write
\begin{align}
\label{eq:w5likeDecomp}
	g = s \Matrix{1&u_2&u_3\\&1&u_1\\&&1}\Matrix{t_1 t_2\\&t_1\\&&1}\Matrix{1\\0&1\\0&z_1&1}w_5\Matrix{1&0&v_3\\&1&v_1\\&&1},
\end{align}
where
\begin{align*}
	s =& B_1(g), & t_1 =& \frac{A_2(g)}{B_1(g)^2}, & t_2 =& \frac{g_{32} \det(g)}{A_2(g)^2}, & z_1 =& \frac{A_1(g)}{B_1(g)}, \\
	u_1 =& \frac{g_{22}}{B_1(g)}, & u_2 =& \frac{g_{11}g_{32}-g_{12}g_{31}}{A_2(g)}, & u_3 =& \frac{g_{12}}{B_1(g)}, & v_1 =& -\frac{B_2(g)}{A_2(g)}, & v_3 =& -\frac{C_2(g)}{A_2(g)}.
\end{align*}

\subsection{The test function}
\label{sect:TestFun}
For $f:\R^\times\to\C$ smooth and compactly supported on $[-T_2, -T_1]\cup[T_1, T_2]$, define $F(g)$ to be zero if $B_1(g) A_2(g) = 0$ and otherwise
\[ F(g) = \psi_I(uv) f\paren{C^3 t_2\sqrt{\abs{t_1}}} h_1(t_1) h_2\paren{X \frac{z_1}{\sqrt{\abs{t_1}}}} h_3(v_1) h_3(v_3), \]
using the coordinates of $g$ as in \eqref{eq:w5likeDecomp}, where each $h_i$ is smooth, $h_1,h_3$ are compactly supported on $(\frac{1}{2},2)$, $h_2$ is even and compactly supported on $(-\frac{1}{200},\frac{1}{200})$, $h_1(t_1)=1$ for $t_1 \in (1-\frac{1}{10},1+\frac{1}{10})$, $h_2(0)=1$, $\int_0^\infty h_3(x) dx = 1$, and $X,C > 0$.
Despite the piece-wise definition, we will show such an $F$ is smooth momentarily.
In \cref{thm:MainThm}, we take $C=X=1$.

Note that there are many other choices one could make, but these generally complicate the analysis and it's not clear that any improvement may be made in the bounds.
In particular, we might consider $h_1$ instead of the form $h_1(X'(t_1-1))$, but the integral in $\wtilde{h_1}(\varepsilon \sqrt{\abs{y_1}})$ of \cref{thm:MainThm} already forces $x_2=1+\BigO{X^{-1}}$ (since $\abs{y_1} \gg X^2$ inside $T_{w_l}(F)$), so we would only save if $X'$ were much larger than $X$, but this loses quite badly in the Weyl law. 

We will now use this test function throughout the rest of the paper, and we collect two facts about $F$ here:
\begin{prop}
\label{prop:Fprops}
For the function $F$ constructed above,
\begin{enumerate}
\item $F$ is smooth and compactly supported in the Iwasawa coordinates on $U(\R)\backslash G$.
\item For $i=1,2$, the action of the Casimir operator $\Delta_i$ on $F$ satisfies $\Delta_i^j F \ll X^j$.
\end{enumerate}
\end{prop}
We defer the proof to \cref{sect:TestFunAnalysis}.

For $y \in Y_{w_5}$, we have
\[ T_{w_5}(F)(y) = f(y_2) h_1(1) h_2(0) \int_{-\infty}^\infty \int_{-\infty}^\infty h_3(v_1) h_3(v_3) dv_1 \, dv_3 = f(y_2). \]
Thus \cref{thm:w5BesselExpand} is an immediate consequence of \cref{prop:Density}.

By the support of $h_1(t_1)$, $T_{w_l}(F)(y)$ is zero unless $\sgn(y_1)=-1$, in which case we have
\begin{align*}
	T_{w_l}(F)(y) =& \int_{U(\R)} \e{\frac{y_1}{x_2}-x_2} f\paren{-\sgn(x_2)C^3 y_2 \sqrt{\abs{y_1}}} h_1\paren{-\frac{y_1}{x_2^2}} h_2\paren{\frac{X}{\sqrt{\abs{y_1}}}} \\
	& \qquad h_3(x_1) h_3(x_3-x_1 x_2) dx \\
	=& \sum_{\varepsilon \in \set{\pm1}} f(\varepsilon C^3 y_2 \sqrt{\abs{y_1}}) h_2\paren{\frac{X}{\sqrt{\abs{y_1}}}} \int_0^\infty \e{-\varepsilon\paren{\frac{y_1}{x_2}-x_2}} h_1\paren{-\frac{y_1}{x_2^2}} dx_2,
\end{align*}
and this is again zero unless $C^3 \sqrt{\abs{y_1}} \abs{y_2} \in (T_1, T_2)$, $\abs{y_1} \ge (200 X)^2$.
Substituting $x_2 \mapsto \sqrt{\abs{y_1}}/x_2$ gives the form of $T_{w_l}(F)(y)$ in \cref{thm:MainThm}.

\subsection{The hyper-Kloosterman sum Kuznetsov formula}
For our test function $F$, we have $T_w(F)=0$ unless $w=w_5,w_l$, so \cref{prop:SpectralInterpret} reduces to the form of \cref{thm:MainThm}, but also the support of $T_{w_l}(F)$ implies the moduli in the long-element term satisfy
\begin{align*}
	\abs{c_2} \le& T_1^{-2/3} C^2 \abs{m_1 m_2^2 n_1^2 n_2}^{1/3}, \\
	1 \le \abs{c_1} \le& \frac{\abs{m_1 n_2 c_2}^{1/2}}{200 X} \le \frac{C \abs{m_1^2 m_2 n_1 n_2^2}^{1/3}}{200 T_1^{1/3} X}.
\end{align*}
Thus for $X$ satisfying \eqref{eq:XBd}, we have $T_{w_l}(F)\paren{mcw_ln^{-1}w_l}=0$ as well and \cref{thm:MainCor} holds.

Again, for this $F$, $T_{w_l}(F)(y)$ doesn't decay rapidly as $y_1 \to -\infty$, so we cannot apply the (long-element) arithmetic Kuznetsov formula (which would cancel the spectral side making the sum of hyper-Kloosterman sums equal to zero).

\subsubsection{The trivial bound on the sum of hyper-Kloosterman sums}
\label{sect:LarsenBd}
Using \cite[App, Thm 1]{BFG}
\[ S_{w_5}(m,n,c) \ll_{m,n,\epsilon} c_1^{2+\epsilon}, \]
the trivial bound on the $w_5$ term is
\begin{align*}
	\sum_{\substack{0 \ne c_1, c_2 \in \Z \\ m_1 c_2=n_2 c_1^2}} S_{w_5}(m,n,c) f\paren{C^3\frac{m_1^2 m_2 n_1}{n_2 c_1^3}} &\ll_{f,m,n,\epsilon} C^{3+\epsilon}.
\end{align*}

\subsubsection{Aside on the long-element term}
\label{sect:LEAside}
In \cref{sect:tildeh1Stationary}, we will show
\begin{align}
\label{eq:TwlFUBd}
	T_{w_l}(F)(y) \ll \abs{y_1}^{1/4}.
\end{align}

For $X$ smaller than in \eqref{eq:XBd}, the trivial bound on the long element term in \cref{thm:MainThm} is
\begin{align*}
	& \sum_{0 \ne c_1, c_2 \in \Z} S_{w_l}(m,n,c) T_{w_l}(F)\paren{\frac{m_1 n_2 c_2}{c_1^2}, \frac{m_2 n_1 c_1}{c_2^2}} \\
	&\ll_{f,m,n,\epsilon} \sum_{\abs{c_2} \asymp C^2} \sum_{\abs{c_1} \ll \sqrt{c_2}/X}\sqrt{(c_1,c_2)} (c_1 c_2)^{\frac{1}{2}+\epsilon} \paren{\frac{c_2}{c_1^2}}^{1/4} \\
	&\ll X^{\epsilon-1} C^{\frac{9}{2}+\epsilon},
\end{align*}
using \cite[(2.10)]{LargeSieve}.
Thus if one could handle a moderately large sum of long-element Kloosterman sums, it is possible to shorten the spectral average.

Note that at $X \asymp_{m,n} C$, the trivial bound on the long-element term matches the trivial bound on the spectral expansion.
Thus it makes sense that we shouldn't be able to improve the trivial bound on the spectral side:
We are encountering the usual uncertainty principle that the trivial bound on the spectral side times the trivial bound on the long-element term is essentially constant, except that the spectral side does not ``see'' that the long-element term disappears for $X$ much larger than $C$.
It is somewhat amusing that in our search for cancellation between terms of the sums of hyper-Kloosterman sums, we've apparently come across cancellation between terms in the spectral sum.

\subsubsection{The spectral side}
\label{sect:SpectralSide}
In this section, we prove \cref{thm:MainInDepth}.

For $d=0,1$, define
\begin{align*}
	&\wtilde{F}^d(\mu) = \\
	& -\frac{\pi^2 C^3}{X} \cscmu^d(\mu) \sum_{\varepsilon_2 \in \set{\pm 1}} \sum_{w \in \Weyl} \sgn(w) \paren{\frac{C}{2\pi}}^{-3\mu^w_1} \chi_d^w((-1,\varepsilon_2)) \frac{\what{f}(\sgn(\Im(\mu^w_2-\mu^w_3)) \varepsilon_2,-1+\mu^w_1)}{\Gamma(1+\mu^w_1-\mu^w_2) \Gamma(1+\mu^w_1-\mu^w_3)} \\
	& \qquad e^{\frac{\pi}{2}\abs{\Im(\mu^w_2-\mu^w_3)}} \what{h}_2\paren{\sgn(\Im(\mu^w_2-\mu^w_3)) \frac{\paren{\mu^w_2-\mu^w_3}^2}{16\pi^2 X}},
\end{align*}
and for $d \ge 2$,
\begin{align*}
	\wtilde{F}^d(\mu^d(r)) =& \frac{i^d \pi C^3}{2 X} \sum_{\varepsilon \in \set{\pm 1}^2} \paren{\frac{C}{2\pi}}^{6r} Q(d,-r) \what{f}(\varepsilon_1 \varepsilon_2,-1-2r) (\varepsilon_1 \varepsilon_2)^d \what{h}_2\paren{\varepsilon_1 \frac{(d-1)^2}{16\pi^2 X}}.
\end{align*}

\begin{prop}
\label{prop:hatFBounds}
Suppose $C,X \ge 1$ and either $\mu$ tempered, i.e. $\mu \in \frak{a}^d_0$, or $\mu$ complementary, i.e. $d=0,1$ and $\mu=(x+it,-x+it,-2it)$, $0<x<\theta$.
We have the trivial bound
\[ \what{F}^d(\mu) \ll (C+X+\norm{\mu})^\epsilon C^3 \times \piecewise{X^{-3/4} & \If \mu \text{ is tempered}, \\ \norm{\mu}^{1+x} C^{3x} X^{2x-1} & \If \mu \text{ is complementary},} \]
and $\what{F}^d(\mu) \ll_\epsilon C^5 \norm{\mu}^{-100}$ unless $\norm{\mu} \ll X^{1/2+\epsilon}$ and  $\min_i \abs{\mu_i} \ll X^{1+\epsilon} \norm{\mu}^{-2}$.

If $1 \ll \norm{\mu} \ll X^{\frac{1}{2}+\epsilon}$, we write
\[ \what{F}^d(\mu) = \what{F}^d_A(\mu)+\what{F}^d_B(\mu)+\what{F}^d_C(\mu)+\what{F}^d_D(\mu), \]
and these satisfy
\begin{enumerate}
\item $\displaystyle \what{F}^d_D(\mu) \ll \piecewise{C^{3+\epsilon} X^{\epsilon-\frac{23}{4}} & \If \mu \text{ is tempered}, \\ C^{3+3\theta+\epsilon} X^{\frac{5}{2}\theta-\frac{11}{2}+\epsilon} & \If \mu \text{ is complementary},}$
\item $\displaystyle \what{F}^d_C(\mu) \ll C^{\frac{3}{2}+\epsilon} X^{-1+\epsilon} \norm{\mu}^{-2}$,
\item $\what{F}^d_A(\mu)+\what{F}^d_B(\mu) \ll X^{-100}$ unless one of $\abs{\Im(\mu_i)} \ll X^\epsilon$, in which case
\[ \what{F}^d_B(\mu) \ll \piecewise{C^{3+\epsilon} X^{\epsilon-2} \norm{\mu}^{-1} & \If \mu \text{ is tempered}, \\ C^{3+3\theta+\epsilon} X^{\epsilon-1} & \If \mu \text{ is complementary},} \]
\item $\what{F}^d_A(\mu)$ is supported on tempered $\mu$ and we have $\what{F}^d_A(\mu) = \wtilde{F}^d(\mu) + \BigO{C^{3+\epsilon} X^{\epsilon-2} \norm{\mu}^{-1}}$, unless $d=0,1$ and some $\abs{\mu_i-\mu_j} < 1$, $i \ne j$, in which case $\what{F}^d_A(\mu) \ll C^{3+\epsilon} X^{\epsilon-1} \norm{\mu}^{-1}$.
\end{enumerate}
\end{prop}
The proof is deferred to \cref{sect:hatFBounds}.
Note that in the trivial bound, if $\mu$ is complementary and $\norm{\mu} < \epsilon$, then also $x < \epsilon$ so the tempered bound dominates.

Now we consider
\[ \abs{\sum_{d=0}^\infty \int_{\mathcal{B}_3^{d*}} \frac{\lambda_\Xi(n)\wbar{\lambda_\Xi(m)}}{L(1,\AdSq\Xi)} \paren{\what{F}^d(\mu_\Xi)-\what{F}^d_A(\mu_\Xi)} d_H\Xi} \ll_{m,n} \sum_{d=0}^\infty \int_{\mathcal{B}_3^{d*}} \frac{\abs{\what{F}^d(\mu_\Xi)-\what{F}^d_A(\mu_\Xi)}}{L(1,\AdSq\Xi)}  d_H\Xi. \]
For the moment we treat the entirety of the cuspidal and continuous spectra by realizing that the contributions of the Eisenstein series are not larger than that of the cusp forms (i.e. comparing Weyl laws).

As there are $\BigO{1}$ forms with $\norm{\mu} < \epsilon$, this part of the spectrum contributes $C^{3+\epsilon} X^{\epsilon-3/4}$.

From the Weyl law \eqref{eq:WeylLawCompSer}, the contribution of the complementary spectrum is bounded by
\[ C^{\frac{3}{2}+\epsilon} X^{-1/2+\epsilon}+ C^{3+3\theta+\epsilon} X^{\epsilon-1}, \]
since $\min_i \abs{\Im(\mu_i)} \ll X^\epsilon$ implies $\norm{\mu} \ll X^\epsilon$ for these forms.

From the Weyl laws \eqref{eq:WeylLaw01} and \eqref{eq:WeylLaw2plus}, the contribution of $\what{F}^d_D(\mu)$ and $\what{F}^d_C(\mu)$ on the tempered spectrum is bounded by
\[ C^{3+\epsilon} X^{\epsilon-\frac{13}{4}}+C^{\frac{3}{2}+\epsilon} X^{1/2+\epsilon}. \]

Lastly, the Weyl law \eqref{eq:WeylLaw01} implies for $T \ge 1$ that
\begin{align*}
	\sum_{\substack{\varphi \in \mathcal{B}^{d*}_{3,\text{cusp}}\\ \norm{\mu_\varphi} \le T \\ \min_i \abs{\mu_i} \ll T^\epsilon}} \frac{1}{L(1,\AdSq \varphi)} \ll T^{4+\epsilon},
\end{align*}
and so the contribution of $\what{F}^d_B(\mu)$ on the tempered spectrum is bounded by $C^{3+\epsilon} X^{-1/2+\epsilon}$.

Next we replace $\what{F}^d_A(\mu)$ with $\wtilde{F}^d(\mu)$.
Away from the Weyl-chamber walls (i.e. $\min_{i\ne j} \abs{\mu_i-\mu_j} < 1$), the error term is the same size as $\what{F}^d_B(\mu)$, so this contribution is bounded by $C^{3+\epsilon} X^{-1/2+\epsilon}$.
Notice that $\abs{\mu_1-\mu_2} < 1$ implies $\mu_3$ is within 1 of $-2\mu_1$, so if we also require $\min_i \abs{\mu_i} \ll X^\epsilon$, then actually $\norm{\mu} \ll X^\epsilon$.
Hence the contribution of $\what{F}^d_A(\mu)$ on the spectrum near the Weyl-chamber walls is bounded by $C^{3+\epsilon} X^{\epsilon-1}$.

Finally, the function $\wtilde{F}^d_0(\mu)$ in \cref{thm:MainInDepth} is the primary asymptotic $\wtilde{F}^d(\mu) = \wtilde{F}^d_0(\mu)+\BigO{C^3 X^{\epsilon-1} \norm{\mu}^{-3}}$.
In the $d=0,1$ case, we used the fact that $\chi_d^{ww_3}((-1,\varepsilon_2)) = \chi_d^w((-1,-\varepsilon_2))$ and we applied Stirling's formula in the form
\[ \Gamma(it) = \sqrt{2\pi} \exp\paren{-\frac{\pi}{2}t-i\frac{\pi}{4}-it} \abs{t}^{it-\frac{1}{2}}\paren{1-\frac{i}{12 t}}\paren{1+ \BigO{\frac{1}{t^2}}}. \]
The error in replacing $\wtilde{F}^d(\mu)$ with $\wtilde{F}^d_0(\mu)$ contributes $\BigO{C^3 X^{\epsilon-\frac{1}{2}}}$ to the error in \cref{thm:MainInDepth}.

On its support, $\wtilde{F}^d_0(\mu)$ is bounded by $C^3 X^{-1} \norm{\mu}^{-1}$, so the tempered spectrum in the range $\norm{\mu} \le X^{1/6}$ contributes at most $\BigO{C^3 X^{\epsilon-\frac{1}{2}}}$, and we have discarded that.

The contribution of minimal parabolic Eisenstein series (which only occurs at $d=0$) is at most
\[ \int_R C^3 X^{-1} \norm{\mu}^{-1} \abs{d\mu} \ll C^3 X^{\epsilon-1}, \qquad R := \set{\mu \setdiv \Re(\mu)=0, \min_i \abs{\mu_i} \ll X^\epsilon, \epsilon < \norm{\mu} \ll X^{\frac{1}{2}+\epsilon}}. \]

The number of $SL(2,\Z)$ cusp forms (either holomorphic modular forms or spherical Maass forms) with spectral parameter (equivalently, the square-root of the Laplacian eigenvalue) less than $X^{\frac{1}{2}+\epsilon}$ is bounded by $X^{1+\epsilon}$, hence the contribution of the maximal parabolic Eisenstein series is at most $C^3 X^{-\frac{1}{2}+\epsilon}$.

\section{Analyzing the test function}
\label{sect:TestFunAnalysis}

\subsection{$F$ in the Iwasawa coordinates}
We now prove \cref{prop:Fprops}.1.
In terms of the Iwasawa coordinates on $U(\R)\backslash G$, we have
\begin{align*}
	F(yk) =& \e{-y_1 \csc \theta_2 \paren{\cos \theta_1 \cot \theta_3 - \cos \theta_2 \sin \theta_1} - y_2 \tan \theta_1 - 
 \csc \theta_2 \sin \theta_3 \paren{\cos \theta_2 + \cot \theta_3 \tan \theta_1}} \\
	& \qquad \times f\paren{-C^3 y_2 \csc \theta_2 \sec^2 \theta_1 \sin \theta_3 \sqrt{-y_1 \cos \theta_1 \csc \theta_2 \csc^2 \theta_3}} \\
	& \qquad \times h_1\paren{-y_1 \cos \theta_1 \csc \theta_2 \csc^2 \theta_3} h_2\paren{X \frac{-\cot \theta_3}{\sqrt{-y_1 \cos \theta_1 \csc \theta_2 \csc^2 \theta_3}}} \\
	& \qquad \times h_3\paren{\csc\theta_2 \sin\theta_3 \paren{\cos\theta_2\cot\theta_3 - \tan \theta_1}} \\
	& \qquad \times h_3\paren{-\csc\theta_2 \sin\theta_3 \paren{\cos\theta_2+\cot\theta_3 \tan \theta_1}} 
\end{align*}
Suppose $yk$ is within the support of $F$, then $\abs{\cot\theta_3} < \frac{1}{100}$ (by the support of $h_2$, since $-y_1 \cos \theta_1 \csc \theta_2 \csc^2 \theta_3 \in (\frac{1}{2},2)$ by the support of $h_1$), so $\abs{\sin\theta_3} > \frac{99}{100}$.
If $\abs{\sin\theta_2} < \frac{1}{3}$, then both $\abs{\cos\theta_2\cot\theta_3 - \tan \theta_1}, \abs{\cos\theta_2+\cot\theta_3 \tan \theta_1} < \frac{200}{297}$ (by the support of $h_3$), but the first gives $\abs{\tan\theta_1} < \frac{203}{297}$ and the second gives $\abs{\cos\theta_2} < \frac{203}{297}$, a contradiction.
Now assume $\abs{\sin\theta_2} \ge \frac{1}{3}$ so we have $\abs{\cos\theta_2\cot\theta_3 - \tan \theta_1}, \abs{\cos\theta_2+\cot\theta_3 \tan \theta_1} \in (\frac{1}{6},\frac{200}{99})$.
The first gives $\abs{\tan\theta_1} \in (\frac{47}{300},\frac{201}{99})$ and the second gives $\abs{\cos\theta_2} > \frac{3}{22}$.

Hence all $\sin\theta_i,\cos\theta_i \asymp 1$ except $\cos\theta_3$ (which may be arbitrarily small), and by the support of $f$ and $h_1$, we also have $y_1\asymp 1, y_2 \asymp C^{-3}$ and $F$ is compactly supported on $U(\R)\backslash G$ in the Iwasawa coordinates.
Also, the only way $F$ could fail to be smooth in the Iwasawa coordinates is if $\cos\theta_3$ occured in a denominator somewhere in our expression for $F(yk)$.
It doesn't.
One might also worry that $\cos\theta_3$ occured in a denominator somewhere in the action of the Lie algebra, see \cite[(31),(39),(73),(74)]{HWII}.
It doesn't.

\subsection{The action of the Lie algebra in the Bruhat-like coordinates}
\label{sect:BruhatLikeLie}
We now prove \cref{prop:Fprops}.2.
The (normalized) Casimir operators are defined by
\begin{align}
\label{eq:CasimirDef}
	\Delta_1 = -\frac{1}{2}\sum_{i,j} E_{ij} \circ E_{ji}, \qquad \Delta_2 = \frac{1}{3}\sum_{i,j,k} E_{ij} \circ E_{jk}\circ E_{ki}+\Delta_1,
\end{align}
where $E_{i,j}$, $i,j\in\set{1,2,3}$ is the element of the Lie algebra of $GL(3,\R)$ whose matrix has a 1 at position $i,j$ and zeros elsewhere, using the standard indexing.

As differential operators acting on smooth functions of $G$, the $E_{ij}$ have the form
\[ E_{11} = t_1 \partial_{t_1}-t_2 \partial_{t_2}-v_3 \partial_{v_3}+z_1 \partial_{z_1}, \qquad E_{21} = -v_3 \partial_{v_1}+\partial_{z_1}, \]
\[ E_{31} = t_1 v_3 \partial_{t_1}-2 t_2 v_3 \partial_{t_2}+t_2 \partial_{u_2}-v_1 v_3 \partial_{v_1}-v_3^2 \partial_{v_3}+(v_1+v_3 z_1) \partial_{z_1}, \]
\[ E_{12} = s z_1 \partial_{s}-2 t_1 z_1 \partial_{t_1}+t_2 z_1 \partial_{t_2}+t_1 \partial_{u_1}+t_1 u_2 \partial_{u_3}-v_1 \partial_{v_3}-z_1^2 \partial_{z_1}, \]
\[ E_{22} = s \partial_{s}-t_1 \partial_{t_1}-v_1 \partial_{v_1}-z_1 \partial_{z_1}, \]
\begin{align*}
	E_{23} =& s (v_1+v_3 z_1) \partial_{s}-t_1 (v_1+2 v_3 z_1) \partial_{t_1}-t_2 (v_1-v_3 z_1) \partial_{t_2}+t_1 v_3 \partial_{u_1}-t_2 z_1 \partial_{u_2} \\
	& \qquad +t_1 (t_2+u_2 v_3) \partial_{u_3}-v_1^2 \partial_{v_1}-v_1 v_3 \partial_{v_3}-z_1 (v_1+v_3 z_1) \partial_{z_1},
\end{align*}
\[ E_{13} = \partial_{v_3}, \qquad E_{23} = \partial_{v_1}, \qquad E_{33} = t_2 \partial_{t_2} + v_1 \partial_{v_1}+v_3\partial_{v_3}. \]

For $d \in \N_0^5$, define
\[ F^{(d)}(g) = \psi_I(uv) f^{(d_1)}\paren{C^3 t_2\sqrt{t_1}} h_1^{(d_2)}(t_1) h_2^{(d_3)}\paren{X \frac{z_1}{\sqrt{t_1}}} h_3^{(d_4)}(v_1) h_3^{(d_5)}(v_3) \]
and let $\set{e_1,\ldots,e_5}$ be the standard basis of $\R^5$, then we have
\begin{align*}
	E_{11}F^{(d)} =& -\frac{1}{2} C^3 \sqrt{t_1} t_2 F^{(d+e_1)}+t_1 F^{(d+e_2)}+\frac{X z_1}{2 \sqrt{t_1}} F^{(d+e_3)}-v_3 F^{(d+e_5)}, \\
	E_{12}F^{(d)} =& -2 t_1 z_1 F^{(d+e_2)}-v_1 F^{(d+e_5)}+2 i \pi  t_1 F^{(d)}, \\
	E_{13}F^{(d)} =& F^{(d+e_5)}, \\
	E_{21}F^{(d)} =& \frac{X}{\sqrt{t_1}} F^{(d+e_3)}-v_3 F^{(d+e_4)}-2 i \pi  v_3 F^{(d)}, \\
	E_{22}F^{(d)} =& -\frac{1}{2} C^3 \sqrt{t_1} t_2 F^{(d+e_1)}-t_1 F^{(d+e_2)}-\frac{X z_1}{2 \sqrt{t_1}} F^{(d+e_3)}-v_1 F^{(d+e_4)}-2 i \pi  v_1 F^{(d)}, \\
	E_{23}F^{(d)} =& F^{(d+e_4)}+2 i \pi  F^{(d)}, \\
	E_{31}F^{(d)} =& -\frac{1}{2} 3 C^3 \sqrt{t_1} t_2 v_3 F^{(d+e_1)}+t_1 v_3 F^{(d+e_2)}+\frac{X}{\sqrt{t_1}}\paren{2v_1+v_3 z_1} F^{(d+e_3)} \\
	& \qquad -v_1 v_3 F^{(d+e_4)}-v_3^2 F^{(d+e_5)}+2 i \pi  F^{(d)} (t_2-v_1 v_3), \\
	E_{32}F^{(d)} =& -\frac{1}{2} 3 C^3 \sqrt{t_1} t_2 v_1 F^{(d+e_1)}-t_1 \paren{v_1+2v_3 z_1} F^{(d+e_2)}-\frac{v_1 X z_1}{2 \sqrt{t_1}} F^{(d+e_3)} \\
	& \qquad -v_1^2 F^{(d+e_4)}-v_1 v_3 F^{(d+e_5)}+2 i \pi \paren{t_1 v_3-t_2 z_1-v_1^2} F^{(d)}, \\
	E_{33}F^{(d)} =& C^3 \sqrt{t_1} t_2 F^{(d+e_1)}+v_1 F^{(d+e_4)}+v_3 F^{(d+e_5)}+2 i \pi  v_1 F^{(d)}.
\end{align*}

So all $E_{ij} F^{(d)}(g) \ll 1$ except $E_{21} F^{(d)}(g), E_{31} F^{(d)}(g) \ll X$, and the same when applying $E_{ij}$ to any of the coefficients.
As neither $E_{21}$ nor $E_{31}$ can occur more than once in the same term of \eqref{eq:CasimirDef} (nor can both occur in the same term), we have the result by applying the symmetries of the Casimir operators as in \cref{prop:Density}.
It might be tempting to say $\norm{\mu} \ll X^{1/3+\epsilon}$, but this is only true when all $\abs{\mu_i} \gg \norm{\mu}$, which is not the important case here.

\subsection{\texorpdfstring{Stationary phase for $\wtilde{h_1}(u)$}{Stationary phase for tilde-h1(u)}}
\label{sect:tildeh1Stationary}
Using \cref{prop:BKY8.2} with $\omega=h_1$, $\phi(t)=u(t+t^{-1})$, $t_0=1$, $T=U_1=U=Q=1$, $Y=\abs{u}$, $\delta=\frac{1}{3}$ and $\phi''(t_0)=2u$ gives
\[ \wtilde{h_1}(u) = \abs{2u}^{-1/2} \e{2u} \sum_{j=0}^{9A} p_j(u) +\BigO{\abs{u}^{-A}}, \]
where
\[ p_j(u) = \frac{\e{\frac{1}{8}\sgn(u)}}{j!} \left. \paren{\frac{i}{8\pi u} \frac{d^2}{dt^2}}^j h_1\paren{t^2} \e{-\frac{u}{t}(t-1)^3} \right|_{t=1}. \]
Since $h_1(1)=1$ and all of the derivatives $h_1^{(j)}(1)=0$ for $j \ge 1$, we can replace $h_1(t^2)$ with 1 in this expression.

One can check that
\[ p_j(u) = \e{\frac{\sgn(u)}{8}} \sum_{k=\floor{(j+2)/3}}^j c_{j,k,\sgn(u)} \abs{u}^{-k} \]
for some constants $c_{j,k,\varepsilon}$, hence
\begin{align}
\label{eq:tildeh1expand}
	\wtilde{h_1}(u) = \abs{2u}^{-1/2} \e{2u+\frac{\sgn(u)}{8}} \sum_{j=0}^N c_{j,\sgn(u)} \abs{u}^{-j} +\BigO{\abs{u}^{-N-\frac{3}{2}}},
\end{align}
where
\[ c_{j,\varepsilon} = \sum_{k=j}^{3j} c_{k,j,\varepsilon}, \]
and in particular, $c_{0,\varepsilon} = 1$.

For $T_{w_l}(F)$, this gives
\begin{align*}
	T_{w_l}(F)(y) =& 2^{-1/2} \sum_{\varepsilon \in \set{\pm1}} \e{2\varepsilon\sqrt{\abs{y_1}}+\frac{\varepsilon}{8}} f(\varepsilon C^3 y_2 \sqrt{\abs{y_1}}) \sum_{j=0}^N \frac{c_{j,\varepsilon}}{\abs{y_1}^{(2j-1)/4}} h_2\paren{\frac{X}{\sqrt{\abs{y_1}}}} \\
	& \qquad + \BigO{\abs{y_1}^{-(2N+1)/4}},
\end{align*}
and in particular, \eqref{eq:TwlFUBd} holds.

\subsection{\texorpdfstring{Bounds for $\what{F}^d(\mu)$}{Bounds for hat-Fd(mu)}}
\label{sect:hatFBounds}
We now proceed to prove \cref{prop:hatFBounds}.

\subsubsection{\texorpdfstring{A direct integral for $\what{F}^d(\mu)$}{A direct integral for hat-Fd(mu)}}
We start with
\begin{align*}
	\what{F}^d(\mu) =& 4 \pi^4 \int_{-\infty}^\infty \int_{-\infty}^\infty T_{w_l}(F)(y) K^{d*}_{w_l}(y,\wbar{\mu}) \frac{dy_1 \, dy_2}{(y_1 y_2)^2},
\end{align*}
where
\[ \abs{16\pi^4 y_1 y_2} K^{d*}_{w_l}(y,\mu) := K^d_{w_l}(y,\mu). \]

For $\mu \in \frak{a}^d_0$, the error in the stationary phase for $\wtilde{h_1}(u)$ and the bound \cite[Lem 22]{ArithKuzII} give
\begin{align*}
	(C+X+\norm{\mu})^\epsilon \int_{\abs{y_1} \gg X^2} \int_{y_2 \asymp C^{-3} \abs{y_1}^{-1/2}} \abs{y_1}^{-(2N+1)/4} \frac{\norm{\mu}^{1/2}}{\abs{y_1}^{1/4}} \frac{dy_1 \, dy_2}{(y_1 y_2)^2} \ll& (C+X+\norm{\mu})^\epsilon C^3 X^{-N-\frac{7}{4}}.
\end{align*}
For $d=0,1$ and $\mu=(x+it,-x+it,-2it)$, $0<x\le\theta$, the bound \cite[Lem 23]{ArithKuzII} gives
\begin{align*}
	& (C+X+\norm{\mu})^\epsilon \norm{\mu}^{1+x} \int_{\abs{y_1} \gg X^2} \int_{y_2 \asymp C^{-3} \abs{y_1}^{-1/2}} \abs{y_1}^{-(2N+1)/4} \abs{y_1}^{\frac{1}{2} (x-\frac{1}{2})} \abs{y_2}^{-x} \frac{dy_1 \, dy_2}{(y_1 y_2)^2} \\
	&\ll (C+X+\norm{\mu})^\epsilon \norm{\mu}^{1+x} C^{3+3x} X^{2x-N-2}.
\end{align*}
It is sufficient for our purposes to take $N=3$, as the Weyl law will contribute $X^{5/2}$; instead, we choose $N=4$ as the above bound is slightly weaker than what we will achieve below.

Taking $N=-1$ in the above analysis corresponds to using the main ($j=0$) term in the expansion \eqref{eq:tildeh1expand} of $\wtilde{h}_1(u)$ and gives the trivial bound stated in the proposition.
The bounds on $\norm{\mu}$ relative to $X$ in the proposition follow by applying the above analysis with $\Delta_i^j F$ in place of $F$ and applying \cref{prop:Fprops}.

\textit{\textbf{We now assume $1 \ll \norm{\mu} \ll X^{\frac{1}{2}+\epsilon}$.}}

Substituting on $y$, we have
\begin{align}
\label{eq:whatFdtruncate}
	\what{F}^d(\mu) =& 2^{5/2} \pi^4 \sum_{\varepsilon \in \set{\pm 1}^2} \sum_{j=0}^4 c_{j,\varepsilon_1} \what{F}^d_j(\varepsilon, \wbar{\mu}) + \what{F}^d_D(\mu),
\end{align}
where
\begin{align*}
	& \what{F}^d_j(\varepsilon, \mu) := \\
	& \frac{C^3}{X^{\frac{1}{2}+j}} \int_0^\infty \int_0^\infty f(\varepsilon_1 \varepsilon_2 y_2) \e{2\varepsilon_1 X y_1+\tfrac{\varepsilon_1}{8}} h_2\paren{y_1^{-1}} K^{d*}_{w_l}\paren{\paren{-X^2 y_1^2,\varepsilon_2 \frac{y_2}{C^3 X y_1}},\mu} \frac{dy_1 \, dy_2}{y_1^{\frac{3}{2}+j} y_2^2},
\end{align*}
and the bound on $\what{F}^d_D(\mu)$ in the proposition follows from the analysis above, realizing that, after applying $\norm{\mu} \ll X^{\frac{1}{2}+\epsilon}$ in the bound for the complementary $\mu$.

\subsubsection{\texorpdfstring{A Mellin-Barnes integral for $\what{F}^d_j(\varepsilon, \mu)$}{A Mellin-Barnes integral for hat-Fdj(epsilon,mu)}}
\label{sect:MBforFdj}
Formally, we have
\begin{align*}
	\what{F}^d_j(\varepsilon, \mu) =& \frac{C^3}{X^{\frac{1}{2}+j}} \int_{\Re(s)=(\frac{1}{2},\frac{1}{2})} (2\pi X)^{-2s_1+s_2} \paren{\frac{C}{2\pi}}^{3s_2} \what{K}^d_{w_l}((-1,\varepsilon_2), s,\mu) \\
	& \qquad \int_0^\infty \int_0^\infty f(\varepsilon_1 \varepsilon_2 y_2) \e{2\varepsilon_1 X y_1} h_2\paren{y_1^{-1}} y_1^{-2s_1+s_2-j-\frac{1}{2}} y_2^{-1-s_2} \frac{dy_1 \, dy_2}{y_1 y_2} \frac{ds}{(2\pi i)^2} \\
	 =& \frac{C^3}{X^{\frac{1}{2}+j}} \int_{\Re(s)=(\frac{1}{2},\frac{1}{2})} (2\pi X)^{-2s_1+s_2} \paren{\frac{C}{2\pi}}^{3s_2} \what{K}^d_{w_l}((-1,\varepsilon_2), s,\mu) \\
	& \qquad \what{f}(\varepsilon_1 \varepsilon_2,-1-s_2) \wtilde{h_2}(\varepsilon_1, \tfrac{1}{2}+j+2s_1-s_2) \frac{ds}{(2\pi i)^2}, \\
	\wtilde{h_2}(\varepsilon, s) :=& \int_0^\infty \e{2\varepsilon X y_1+\tfrac{\varepsilon}{8}} h_2\paren{y_1^{-1}} y_1^{-s} \frac{dy_1}{y_1}.
\end{align*}
The formally computation is justified by the Plancherel identity for the Mellin transform; we do not need to worry about bending contours slightly to the left of the $\Re(s_i)=0$ lines since we have some decay coming from $\what{f}$ and $\wtilde{h_2}$.

\subsubsection{\texorpdfstring{Bounds for $\wtilde{h_2}(\varepsilon, s)$}{Bounds for tilde-h2(epsilon,s)}}
Repeated integration by parts shows that $\wtilde{h_2}(\varepsilon, s)$ is entire in $s$ and $\BigO{X^{-100}}$ unless $\abs{s} \ge X^{1-\epsilon}$ since
\[ \frac{d^k}{dy_1^k} h_2\paren{y_1^{-1}} y_1^{-s-1} \ll_k y_1^{-\Re(s)-1} \paren{\frac{1+\abs{s}}{y_1}}^k. \]
Now let $\sigma+it=s$ with $\sigma$ in some fixed compact subset of $\R$ (say $[-10,10]$) and $t$ real with $\abs{t} \ge X^{1-\epsilon}$.

We cut the integral into two pieces with a partition of unity $h_4\paren{y_1/\frac{\abs{t}^{1+\epsilon}}{X}}$ (and its complement $1-h_4$) where $h_4$ is supported on $[0,2]$ and is one on $[0,\frac{3}{2}]$.
Again, repeated integration by parts shows the integral is $\BigO{(X+\abs{s})^{-100}}$ on the region $y_1 \gg \frac{\abs{t}^{1+\epsilon}}{X}$.

In case $\frac{\abs{t}^{1-\epsilon}}{X} > 100$, we again cut the remaining integral into two pieces with a partition of unity $h_4\paren{y_1/\frac{\abs{t}^{1-\epsilon}}{X}}$.
On the region $y_1 \ll \frac{\abs{t}^{1-\epsilon}}{X}$, repeated integration by parts shows the integral is $\BigO{(X+\abs{s})^{-100}}$ since
\[ \paren{y_1 \frac{d}{dy_1}}^k \e{2\varepsilon X y_1} h_2\paren{y_1^{-1}} \ll \abs{t}^{k(1-\epsilon)}. \]

Lastly, we insert a dyadic partition of unity
\[ 1 = \sum_{\ell \in \Z} h_5\paren{y_1 / \frac{2^\ell \abs{t}}{2\pi X}}, \]
where $h_5$ is supported on $[\frac{1}{2},2]$ and one on $(1-\frac{1}{10},1+\frac{1}{10})$.
We may truncate the partition to $\abs{\ell} \le \epsilon \log_2 X$ since otherwise, we are outside the support of
\[ h_4\paren{\frac{X y_1}{\abs{t}^{1+\epsilon}}} \paren{1-h_4\paren{\frac{X y_1}{\abs{t}^{1-\epsilon}}}}. \]
Note: \cref{prop:BKY8.2} does not allow an epsilon-power difference between the $Y$ and $Q$ in the second derivative and the $Y$ and $Q$ in the higher derivatives, hence the need for the dyadic partition of unity.

So we are left with the oscillatory integral
\[ \wtilde{h_2}(\varepsilon, s) = \sum_{\abs{\ell} \le \epsilon \log_2 X} \int_0^\infty \e{\phi(y_1)} \omega_\ell(y_1) dy_1+\BigO{(X+\abs{s})^{-100}} \]
where
\begin{align*}
	\phi(y_1) =& 2\varepsilon X y_1-\frac{t}{2\pi} \log y_1, \\
	\omega_\ell(y_1) =& h_2\paren{y_1^{-1}} y_1^{-\sigma-1} h_4\paren{\frac{X y_1}{\abs{t}^{1+\epsilon}}} \paren{1-h_4\paren{\frac{X y_1}{\abs{t}^{1-\epsilon}}}} h_5\paren{y_1 / \frac{2^\ell \abs{t}}{2\pi X}},
\end{align*}
and the derivatives of the phase and weight functions satisfy
\begin{align*}
	\omega_\ell^{(k)}(y_1) \ll_k& \paren{\frac{2^\ell \abs{t}}{X}}^{-\sigma-1} \paren{\frac{\abs{t}^{1-\epsilon}}{X}}^{-k}, \\
	\phi'(y_1) =& 2\varepsilon X-\frac{t}{2\pi y_1}, \\
	\abs{\phi''(y_1)} \asymp& \frac{\abs{t}}{y_1^2} \asymp \abs{t} \paren{\frac{2^\ell \abs{t}}{X}}^{-2}, \\
	\phi^{(k)}(y_1) \asymp& \frac{\abs{t}}{y_1^k} \asymp \abs{t}\paren{\frac{2^\ell \abs{t}}{X}}^{-k}, \qquad k \ge 2.
\end{align*}
(In case $\frac{\abs{t}^{1-\epsilon}}{X} \le 100$, the factor $1-h_4$ is irrelevant.)

If $\sgn(t) \ne \varepsilon$ or $\ell \ne 0$, the oscillatory integral is $\BigO{(X+\abs{s})^{-100}}$ by \cref{lem:BKY8.1} using $T = \paren{\frac{2^\ell \abs{t}}{X}}^{-\sigma-1}$, $U = \frac{\abs{t}^{1-\epsilon}}{X}$, $R=X$, $Y=\abs{t}$, $Q=\frac{\abs{t}^{1-\epsilon}}{X}$.
Otherwise, the oscillatory integral has a stationary point at $y_1 = \frac{\abs{t}}{2\pi X}$ and \cref{prop:BKY8.2} with $T = \paren{\frac{\abs{t}}{X}}^{-\sigma-1}$, $U = \frac{\abs{t}^{1-\epsilon}}{X}$, $U_1=2\frac{\abs{t}}{X}$, $Y = \abs{t}$, $Q = \frac{\abs{t}}{X}$, $\delta=\frac{1}{3\Max{1,-\sigma-1}}$ gives
\begin{align}
\label{eq:tildeh2Bd}
	\wtilde{h_2}(\varepsilon, s) \ll& \frac{Q\cdot \paren{\frac{\abs{t}}{X}}^{-\sigma-1}}{\sqrt{Y}} = X^\sigma \abs{t}^{-\sigma-\frac{1}{2}}.
\end{align}
By the same analysis, we have
\begin{align}
\label{eq:tildeh2DervBd}
	\frac{d}{ds} \wtilde{h_2}(\varepsilon, s) \ll X^\sigma \abs{t}^{\epsilon-\sigma-\frac{1}{2}}.
\end{align}

\subsubsection{Contour shifting}
In $\what{F}^d_j(\varepsilon, \mu)$, shift $\Re(s_2)$ back to $-\frac{1}{2}+\epsilon$ (with $0 < \epsilon < \frac{1}{14}$ to avoid any poles) and call the shifted contour $\what{F}^d_{j,C}(\varepsilon, \mu)$, then for $d=0,1$, we have
\begin{align*}
	\what{F}^d_{j,C}(\varepsilon, \mu) \ll& C^{\frac{3}{2}} X^\epsilon \int_{\substack{\Re(s_2) = -\frac{1}{2}+\epsilon\\ \abs{s_2} \ll X^\epsilon}} \int_{\substack{\Re(s_1) = \frac{1}{2}\\ \abs{s_1} \ge X^{1-\epsilon}}} \abs{s_1}^{-2-j} \prod_{i=1}^3 \paren{1+\abs{\mu_i}}^{-1+\Re(\mu_i)} \abs{ds_1}\, \abs{ds_2} \\
	\ll& C^{\frac{3}{2}} X^{-1-j+\epsilon} \prod_{i=1}^3 \paren{1+\abs{\mu_i}}^{-1+\Re(\mu_i)},
\end{align*}
and for $d \ge 2$, we have
\begin{align*}
	\what{F}^d_{j,C}(\varepsilon, \mu) \ll& C^{\frac{3}{2}} X^\epsilon \int_{\substack{\Re(s_2) = -\frac{1}{2}+\epsilon\\ \abs{s_2} \ll X^\epsilon}} \int_{\substack{\Re(s_1) = \frac{1}{2}\\ \abs{s_1} \ge X^{1-\epsilon}}} \abs{s_1}^{-2-j} \paren{1+\abs{r}}^{-1} \paren{d+\abs{r}}^{-2} \abs{ds_1}\, \abs{ds_2} \\
	\ll& C^{\frac{3}{2}} X^{-1-j+\epsilon} \paren{1+\abs{r}}^{-1} \paren{d+\abs{r}}^{-2}.
\end{align*}
Note that for the complementary spectrum we have
\[ \prod_{i=1}^3 \paren{1+\abs{\mu_i}}^{\Re(\mu_i)} \asymp \prod_{i=1}^3 \norm{\mu}^{\Re(\mu_i)} = 1 \]
and otherwise on $d=0,1$, we have $\Re(\mu)=0$.
Finally, sacrificing strength for simplicity, we have
\[ \norm{\mu}^{-2} \gg \piecewise{\displaystyle \prod_{i=1}^3 \paren{1+\abs{\mu_i}}^{-1} & \If d=0,1, \\ \displaystyle \paren{1+\abs{r}}^{-1} \paren{d+\abs{r}}^{-2} & \If d \ge 2.} \]
Thus we define
\begin{align*}
	\what{F}^d_C(\mu) =& 2^{5/2} \pi^4 \sum_{\varepsilon \in \set{\pm 1}^2} \sum_{j=0}^3 c_{j,\varepsilon_1} \what{F}^d_{j,C}(\varepsilon, \wbar{\mu}),
\end{align*}
and the bound in the proposition follows.

For the moment, \textit{\textbf{we make the technical assumption that}}
\[ \min_{i\ne j} \abs{\mu_i - \mu_j} > \epsilon \]
for some $\epsilon > 0$ (note that this holds trivially for $d \ge 2$); this avoids some complications with double poles of $\what{K}^d_{w_l}$ and poles of the residues of $\what{K}^d_{w_l}$.
We will discuss the modifications required to remove this technical assumption in \cref{sect:FormsNearWalls}.
Now let $\mathcal{Z}$ be the (finite) collection of poles in $s_2$ to the right of $\Re(s_2)=\frac{1}{2}$; we have $\what{F}^d_j(\varepsilon, \mu) = \what{F}^d_{j,B}(\varepsilon, \mu) + \what{F}^d_{j,C}(\varepsilon, \mu)$, where
\begin{align*}
	\what{F}^d_{j,B}(\varepsilon, \mu) :=& \frac{C^3}{X^{\frac{1}{2}+j}} \sum_{z \in \mathcal{Z}} \paren{\frac{C^3 X}{4\pi^2}}^z \what{f}(\varepsilon_1 \varepsilon_2,-1-z) \int_{\Re(s_1)=\frac{1}{2}} (2\pi X)^{-2s_1} \\
	& \qquad \wtilde{h_2}(\varepsilon_1, \tfrac{1}{2}+j+2s_1-z) \res_{s_2=z} \what{K}^d_{w_l}((-1,\varepsilon_2), s,\mu) \frac{ds_1}{2\pi i},
\end{align*}
since our assumptions rule out any double poles.

For $d=0,1$, we have $\mathcal{Z}=\set{-\mu^w_1 \setdiv w \in \Weyl_3}$, and for $d \ge 2$, we have $\mathcal{Z}=\set{2r}$.
Clearly $\what{F}^d_{j,B}(\varepsilon, \mu) \ll X^{-100}$ unless $\abs{\Im(z)} \ll X^\epsilon$ for some $z \in \mathcal{Z}$ by the super-polynomial decay of $\what{f}(\varepsilon_1 \varepsilon_2,-1-z)$; this will imply the first bound on $\what{F}^d_A(\mu)+\what{F}^d_B(\mu)$ in the proposition (see below for the choice of these two functions).

Now we reverse the Mellin-Plancherel identity using the known Mellin transform of the $J$-Bessel function \cite[6.422.6]{GradRyzh}.
For $d=0,1$, this gives
\begin{equation}
\label{eq:hatFBint01}
\begin{aligned}
	&\what{F}^d_{j,B}(\varepsilon, \mu) = \\
	& -\frac{C^3 \cscmu^d(\mu)}{2\pi X^{\frac{1}{2}+j}} \sum_{w \in \Weyl} \sgn(w) \paren{\frac{C}{2\pi}}^{-3\mu^w_1} \chi_d^w((-1,\varepsilon_2)) \frac{\what{f}(\varepsilon_1 \varepsilon_2,-1+\mu^w_1) H_j(\varepsilon_1, \mu^w_2-\mu^w_3)}{\Gamma(1+\mu^w_1-\mu^w_2) \Gamma(1+\mu^w_1-\mu^w_3)},
\end{aligned}
\end{equation}
while for $d\ge 2$, this gives
\begin{equation}
\label{eq:hatFBint2plus}
\begin{aligned}
	\what{F}^d_{j,B}(\varepsilon, \mu) =& \frac{\varepsilon_2^d}{4\pi^2} \frac{C^3}{X^{\frac{1}{2}+j}} \paren{\frac{C}{2\pi}}^{6r} Q(d,-r) \what{f}(\varepsilon_1 \varepsilon_2,-1-2r) H_j(\varepsilon_1, d-1),
\end{aligned}
\end{equation}
where
\[ H_j(\varepsilon,\nu) := \int_0^\infty \e{2\varepsilon X y_1+\tfrac{\varepsilon}{8}} h_2\paren{y_1^{-1}} J_\nu(4\pi X y_1) \frac{dy_1}{y_1^{3/2+j}}. \]

\subsubsection{The complementary series}
\label{sect:CompSer}
If $d=0,1$ and $\mu=(x+it,-x+it,-2it)$ with $0<x<\theta$ and $\abs{t} \ll X^\epsilon$, we apply \eqref{eq:JsigmaitBd} so that
\begin{align*}
	\what{F}^d_{j,B}(\varepsilon, \mu) \ll& X^\epsilon \frac{C^{3+3\theta}}{X^{1+j}},
\end{align*}
hence we take 
\begin{align}
\label{eq:hatFBCompSer}
	\what{F}^d_B(\mu) =& 2^{5/2} \pi^4 \sum_{\varepsilon \in \set{\pm 1}^2} \sum_{j=0}^3 c_{j,\varepsilon_1} \what{F}^d_{j,B}(\varepsilon, \wbar{\mu})
\end{align}
and $\what{F}^d_A(\mu) = 0$ in this case.

\subsubsection{The tempered forms}
For $\mu \in \frak{a}^d_0$, we have the trivial bound (recall $\abs{\mu_1^w} \ll X^\epsilon$ in the gamma functions of \eqref{eq:hatFBint01}; similarly $\abs{r} \ll X^\epsilon$ in \eqref{eq:hatFBint2plus})
\begin{align}
\label{eq:hatFBtemperedBd}
	\what{F}^d_{j,B}(\varepsilon, \mu) \ll& X^\epsilon \frac{C^3}{X^{1+j}} \norm{\mu}^{-1},
\end{align}
so in this case, we take $\what{F}_B(\mu)$ to be the $j \ge 1$ terms of \eqref{eq:hatFBCompSer} and $\what{F}^d_A(\mu)$ to be the $j=0$ term.

Now we apply the asymptotic expansions \eqref{eq:JitFirstTerm} or \eqref{eq:Jdm1FirstTerm} to $H_0(\varepsilon,\nu)$ with $\nu=it$ or $\nu = d-1$, respectively.
So we need to analyze
\begin{align*}
	H_0(\varepsilon,\nu) =& \sum_\pm \frac{D_{\nu,\pm}}{\sqrt{8 X} \pi} \int_0^\infty \e{2\varepsilon X y_1+\frac{\varepsilon}{8}\pm\paren{2 X y_1+\frac{\nu^2}{16\pi^2 X y_1}+\frac{1}{8}}} h_2\paren{y_1^{-1}} \frac{dy_1}{y_1^{2}} \\
	& \qquad +\BigO{e^{\frac{\pi}{2}\abs{\Im(\nu)}} X^{\epsilon-3/2}},
\end{align*}
where
\[ D_{it,\pm} = e^{\frac{\pi}{2}\abs{t}} \delta_{\pm=-\sgn(t)}, \qquad D_{d-1,\pm} = (\mp i)^d. \]

We can save a factor of $X^{1-\epsilon}$ by integration by parts (since $\nu^2 \ll X^{1+\epsilon}$) unless $\pm = -\varepsilon$, so we have
\begin{align*}
	H_0(\varepsilon,\nu) =& \frac{D_{\nu,-\varepsilon}}{\sqrt{8 X} \pi} \what{h}_2\paren{\varepsilon \frac{\nu^2}{16\pi^2 X}} +\BigO{e^{\frac{\pi}{2}\abs{\Im(\nu)}} X^{\epsilon-3/2}}.
\end{align*}
In both cases, the error term contributes $\BigO{C^3 X^{\epsilon-2} \norm{\mu}^{-1}}$ to $\what{F}^d_A(\mu)$.

\subsubsection{Forms near the Weyl-chamber walls}
\label{sect:FormsNearWalls}
We now discuss how to remove the technical assumption that $\abs{\mu_i-\mu_j} > \epsilon$ for $i \ne j$.
We leave the case $\norm{\mu} < \epsilon$ to the trivial bound, and we stop at \eqref{eq:hatFBtemperedBd}, as taking $j=0$ there is sufficient for this class of forms.

The fundamental issue here is that the sines in the denominator of $\what{K}^d_{w_l}(s,v,\mu)$ (recall \eqref{eq:Kwl1Sgn}) for $d=0,1$ may be zero or arbitrarily close to it; this also gives rise to double poles in the contour shifting.
On the other hand, the original function $K^d_{w_l}(y,\mu)$ is entire in $\mu$; in particular, it is continuous and moreover differentiable near each $\mu_i=\mu_j$.
In fact, we will see that each individual term $\what{F}^d_{j, \cdot}(\varepsilon,\mu)$ of the previous analysis is continuous (so we can deal with $\mu_i=\mu_j$ by taking limits) and differentiable (so we can deal with $\abs{\mu_i-\mu_j} \le \epsilon$ by applying the Mean Value Theorem) there.
So we apply the Mean Value Theorem in the form
\begin{align}
\label{eq:MVT}
	\sup_{\abs{h} \le \epsilon} \abs{\frac{f(h)}{\sin(h)}+\frac{f(-h)}{\sin(-h)}} \le 2\sup_{\abs{h} \le \epsilon} \abs{f'(h)}
\end{align}
(including the ill-defined case where $h=0$ on the left-hand side by continuity), and write $\mu$ in the form, say, $\mu=a(1,1,-2)$ as $\mu = \lim_{h \to 0} a(1,1,-2)+h(1,-1,0)$.
Note that for $\norm{\mu} \ge \epsilon$, there can be at most one pair $\abs{\mu_i-\mu_j} < \epsilon/3$, so we need not worry about second- or third-order derivatives.

To see the differentiability, note that in \eqref{eq:Kwl1Sgn}, for $d=0$, we have $\chi_d^w(v)=1$ and $G^{v_1,v_2}(s,\mu^w)$ is invariant under $w \mapsto w w'$ for some odd permutation $w'$, hence we can inflate the sum over $w \in \Weyl_3$ to a sum over $w \in \Weyl$ and pair terms corresponding to the transposition $(\mu_i, \mu_j) \stackrel{w}{\to} (\mu_j,\mu_i)$ as in the left-hand side of \eqref{eq:MVT}.
For $d=1$, we only need to check the transposition $w_2$ which interchanges $\mu_1$ with $\mu_2$ as that is the only pole of $\cscmu^d(\mu)$ on (or near) the spectrum; it requires a bit more work, but in each case we can see the necessary invariance, e.g. $G^{+-}(s,\mu) = G^{+-}(s,\mu^{w_3})$ so in
\[ (-1)^d G^{+-}(s,\mu)+(-1)^d G^{+-}(s,\mu^{w_5})+G^{+-}(s,\mu^{w_4}), \]
the last term is invariant under interchanging $\mu_1$ and $\mu_2$, while that operation merely swaps the first two terms.
(These facts are rather more obvious if one starts with the power-series definition; see \cite[(82)-(84)]{ArithKuzII}.)
This behavior necessarily carries through to the residues of $\what{K}^d_{w_l}(s,v,\mu)$ in \cref{prop:hatKBd01,prop:hatKBd2plus} and \eqref{eq:hatFBint01}, \eqref{eq:hatFBint2plus}.

From well-known bounds for the digamma function, e.g. applying \cite[8.361.3]{GradRyzh}, we have
\[ \Gamma'(z) \ll \Gamma(z)\paren{1+\abs{\log z}}, \]
so the bounds of \cref{prop:hatKBd01} and \cref{prop:hatKBd2plus} remain true for the derivatives (in $\mu$ using $h \to 0$ as above), at the cost of an epsilon power of $X$.
We may replace \eqref{eq:tildeh2Bd} with \eqref{eq:tildeh2DervBd}, again at a cost of an epsilon power of $X$.
Finally, we may replace the bounds on $J_\nu(x)$ of \cref{lem:BesselBd} with those of \cref{lem:BesselDervBd} at a cost of an epsilon power of $X$.

\section{Bounds on Bessel functions}
\subsection{Simple global bounds on the direct functions}
\label{sect:BesselBd}
In this section, we prove \cref{lem:BesselBd}.

For $x>0$ and $t\in\R$, \cite[7.13.2 (17)]{EMOT} gives
\begin{align*}
	J_{it}(x) =& \frac{1}{\sqrt{2\pi}} (x^2+t^2)^{-1/4} e^{\frac{\pi}{2}\abs{t}} \exp\paren{i\sgn(t) \omega(x,t)} \\
	& \qquad \times\paren{\sum_{j=0}^N \frac{\tilde{a}_j(1+x^2/t^2)}{(x^2+t^2)^{j/2}}+\BigO{x^{-N-1}}}, \\
	\omega(x,t) :=& \abs{t} \arcsinh\frac{\abs{t}}{x}-\sqrt{t^2+x^2}-\frac{\pi}{4},
\end{align*}
where 
\[ \tilde{a}_0(u) = 1, \qquad \tilde{a}_1(u) = -\frac{i}{24}\paren{3-\frac{5}{u}}, \qquad \ldots \]
satisfies $\tilde{a}_j(1+x^2/t^2) \ll_j 1$.
In particular,
\begin{align}
\label{eq:EMOTJitFirstTerm}
	J_{it}(x) =& \frac{1}{\sqrt{2\pi}} (x^2+t^2)^{-1/4} e^{\frac{\pi}{2}\abs{t}} \exp\paren{i\sgn(t) \omega(x,t)}\paren{1 +\BigO{x^{-1}}},
\end{align}
and for $\abs{t} \ll x^{\frac{1}{2}+\epsilon}$, we have
\[ \omega(x, t) = -x+\frac{t^2}{2x}-\frac{\pi}{4}+\BigO{x^{\epsilon-1}}, \]
so
\begin{align}
\label{eq:JitFirstTerm}
	J_{it}(x) =& \frac{1}{\sqrt{2\pi x}} e^{\frac{\pi}{2}\abs{t}} \exp \paren{i\sgn(t) \paren{-x+\frac{t^2}{2x}-\frac{\pi}{4}}} +\BigO{x^{\epsilon-3/2} e^{\frac{\pi}{2}\abs{t}}}.
\end{align}

The $J$-Bessel function at real order is related to the Hankel functions via $J_\sigma(x) = \Re\paren{H^{(1)}_\sigma(x)}$, so for $x > \sigma > 0$, \cite[7.13.2 (11)]{EMOT} gives
\begin{align}
\label{eq:EMOTJp}
	J_\sigma(x) =& \Re\Biggl(\sqrt{\frac{2}{\pi}} \paren{x^2-\sigma^2}^{-1/4} \exp \paren{i \wtilde{\omega}(x,\sigma)} \paren{\sum_{j=0}^N \frac{\wtilde{a}_j(1-x^2/\sigma^2)}{(x^2-\sigma^2)^{j/2}}+\BigO{x^{-N-1}}}\Biggr),
\end{align}
where
\begin{align*}
	\wtilde{\omega}(x,\sigma) :=& \paren{x^2-\sigma^2}^{1/2}+\sigma\arcsin\paren{\sigma/x}-\tfrac{\pi}{2}(\sigma+\tfrac{1}{2}).
\end{align*}
In particular, for $0 < d \ll x^{\frac{1}{2}+\epsilon}$,
\begin{align}
\label{eq:Jdm1FirstTerm}
	J_{d-1}(x) =& \sqrt{\frac{2}{\pi x}} \cos\paren{x+\frac{(d-1)^2}{2x}-\frac{\pi}{2}d+\tfrac{\pi}{4}}+\BigO{x^{\epsilon-3/2}}.
\end{align}

For the first Hankel function, the asymptotics \cite[7.13.2 (13),(15),(22)]{EMOT} imply
\begin{align}
\label{eq:HankelBd}
	H^{(1)}_\sigma(x) \ll \Min{\sigma^{-1/3}, \abs{x^2-\sigma^2}^{-1/4}}, \qquad H^{(1)}_{it}(x) \ll (x^2+t^2)^{-1/4} e^{\frac{\pi}{2}\abs{t}}, \qquad x \ge 1.
\end{align}
From \eqref{eq:HankelBd} and the power series \cite[8.402]{GradRyzh}, we see that $F(\nu) := (x^2-\nu^2-i)^{1/4} e^{i\frac{\pi}{2} \nu} J_\nu(x) \ll 1$ holds on the boundary of the sector $\set{\nu \in \C \setdiv 0 \le \arg \nu \le \frac{\pi}{2}}$ uniformly on $x > 0$.
For a fixed $x$, the bound $F(\nu) \ll_x \exp \paren{\abs{\nu} \paren{1+\abs{\log (x/2)}}}$ holds inside the sector (by the power series expansion), so by Phragm\'en-Lindel\"of for a sector, we have \eqref{eq:PhragLindJBd}.
Also, using \cite[8.471.1]{GradRyzh}, for $-1 < \Re(\nu) < 0$, we have
\[ J_\nu(x) = \frac{2\nu}{x} J_{\nu+1}(x)-J_{\nu+2}(x), \]
and this implies \eqref{eq:JsigmaitBd}.

We may express the $Y$-Bessel function in terms of the Hankel functions as $i Y_{it}(x) = H^{(1)}_{it}(x)-J_{it}(x)$, so \eqref{eq:HankelBd}, the bounds on $J_{it}(x)$, and the above arguments give \eqref{eq:YsigmaitBd}.
Note that on $0 < x \le 1$, the $Y$-Bessel function displays an additional logarithmic dependence on $x$, so we simply avoid that here.

\subsection{Bounds on the order derivative of the $J$-Bessel function}
\label{sect:BesselDervBd}
In this section, we prove \cref{lem:BesselDervBd}.
First, let us recall some well-known asymptotics for the $Y$-Bessel functions:
From \cite[10.8.1, 10.8.2, 10.17.4]{DLMF}, some simple, global bounds for $Y_0$ and $Y_1$ are
\begin{align}
\label{eq:GlobalY0}
	Y_0(y) \ll& \Min{1+\abs{\log y},y^{-1/2}}, \\
\label{eq:GlobalY1}
	Y_1(y) \ll& y^{-1}+\abs{\log y} y^{-1/2},
\end{align}
and more precisely, we have \cite[10.17.4]{DLMF}
\begin{align}
\label{eq:LargeArgY0}
	Y_0(y) =& \sqrt{\frac{2}{\pi y}} \cos \paren{y-\frac{\pi}{4}} + \BigO{y^{-3/2}}.
\end{align}

As in the previous section, we prove the bound for $\nu \ge 0$, then for $\nu$ purely imaginary and Phragm\'en-Lindel\"of implies the lemma.
More precisely, we first show the bound for $x < \abs{t}^{1/3}$ and then for $\abs{t} \le 2$ to avoid some technical complications in the difficult case of $\nu$ purely imaginary.

When $x < \abs{t}^{1/3}$, we can simply use the power series expansion for (the $\nu$-derivative of) $J_\nu(x)$; one can easily derive the bound $\BigO{\log(3+x^{-1})\paren{1+\abs{t}}^{-1/2}\log(3+\abs{t})}$ in this case.
When $\abs{t} \le 2$, we apply a formula of Dunster \cite[(2.1), (2.8)]{Dunster01}:
\begin{align*}
	\frac{\partial}{\partial\nu} J_\nu(x) =& \frac{\pi}{2} Y_\nu(x) +\pi \nu \paren{J_\nu(x) \int_x^\infty J_\nu(t) Y_\nu(t) \frac{dt}{t}-Y_\nu(x) \int_x^\infty \paren{J_\nu(t)}^2 \frac{dt}{t}},
\end{align*}
for $x \ne 0$, $\abs{\arg(x)} < \pi$, provided the unbounded parts of the contours lie on the positive real axis and the contours avoid $(-\infty,0]$.
This formula is effective in the case $\abs{\Im(\nu)} \ll 1$, but the opposite case requires the recovery of exponential factors resulting from the products of three Bessel functions.
Again, a simple application of the $J$- and $Y$-Bessel bounds gives the lemma in this case, with the integrals adding at most a logarithmic factor in case $x$ is near zero.

Now we handle the case where $\nu$ is purely imaginary.
In \cite[eq. (42)]{ApelKrav}, we have
\begin{align*}
	\frac{\partial}{\partial \nu}J_\nu(x) =& \frac{\pi}{2} \nu \int_0^x Y_0(x-y) J_\nu(y) \frac{dy}{y}.
\end{align*}
This integral representation does not apply at $\Re(\nu)=0$, but if we take $0 < \eta < x$ and $N \ge 0$, we can write this as
\begin{align*}
	\frac{\partial}{\partial \nu}J_\nu(x) =& \frac{\pi}{2}\nu \int_\eta^x Y_0(x-y) J_\nu(y) \frac{dy}{y}+\frac{\pi}{2}\nu \int_0^\eta Y_0(x-y) \paren{J_\nu(y)-\sum_{k=0}^N \frac{(-1)^k (y/2)^{\nu+2k}}{k! \Gamma(k+1+\nu)}} \frac{dy}{y} \\
	& \qquad + \frac{\pi}{2}\nu \sum_{k=0}^N \frac{(-1)^k 2^{-\nu-2k}}{k! \Gamma(k+1+\nu)} \int_0^\eta Y_0(x-y) y^{\nu+2k-1} dy.
\end{align*}
Then integration by parts gives
\begin{align*}
	\frac{\partial}{\partial \nu}J_\nu(x) =& \frac{\pi}{2}\nu \int_\eta^x Y_0(x-y) J_\nu(y) \frac{dy}{y}+\frac{\pi}{2}\nu \int_0^\eta Y_0(x-y) \paren{J_\nu(y)-\sum_{k=0}^N \frac{(-1)^k (y/2)^{\nu+2k}}{k! \Gamma(k+1+\nu)}} \frac{dy}{y} \\
	& + \frac{\pi}{2}\nu \sum_{k=0}^N \frac{(-1)^k (\eta/2)^{\nu+2k}}{k! (\nu+2k) \Gamma(k+1+\nu)} \paren{Y_0(x-\eta)-2\eta \int_0^1 Y_1(x-\eta y) y^{\nu+2k} dy},
\end{align*}
and this has analytic continuation to $\Re(\nu)=0$.

Now suppose $x \ge 4 t^{1/3}, t \ge 3$ and $\nu=it$ and in place of the sharp cut-off, we apply a smooth, dyadic partition of unity as $y$ or $x-y$ becomes small, stopping when $y < t^{1/3}$ or $x-y < x/t$ (more precisely, we apply the partition, which depends on $t$, then analytically continue to $\nu=it$), so that
\begin{align*}
	\frac{2}{\pi it} \left.\frac{\partial}{\partial \nu}J_\nu(x)\right|_{\nu=it} =& T_0+\sum_{t^{1/3} \le C \le x/4} T_{1a}(C)+\sum_{x/t \le C \le x/2} T_{1b}(C)+T_{2a} - T_{2b}, \\
	T_0 :=& \int_0^x Y_0(y) J_{it}(x-y) f_0\paren{t y/x} \frac{dy}{x-y}, \\
	T_{1a}(C) :=& \int_0^x Y_0(x-y) J_{it}(y) f_1(y/C) \frac{dy}{y}, \\
	T_{1b}(C) :=& \int_0^x Y_0(x-y) J_{it}(y) f_1((x-y)/C) \frac{dy}{y}, \\
	T_{2a} :=& \int_0^x Y_0(x-y) \paren{J_{it}(y)-\sum_{k=0}^N \frac{(-1)^k (y/2)^{it+2k}}{k! \Gamma(k+1+it)}} f_0\paren{y/t^{1/3}} \frac{dy}{y}, \\
	T_{2b} :=& \sum_{k=0}^N \frac{(-1)^k \paren{t^{1/3}/2}^{it+2k}}{k! (it+2k) \Gamma(k+1+it)} \\
	& \qquad \times \int_0^1 \paren{t^{1/3} Y_1\paren{x-t^{1/3} y} f_0(y)+Y_0\paren{x-t^{1/3}y} f_0'(y)} y^{it+2k} dy,
\end{align*}
where each $f_i$ is smooth, $f_i^{(n)}(y) \ll_n 1$ and, say, $f_0$ supported on $[-1,1]$ (in constructing the partition, $f_0$ will be nonzero at 0) while $f_1$ supported on $[\frac{1}{2},\frac{3}{2}]$; the $C$ sums are dyadic.

For the first term, we use \eqref{eq:GlobalY0} so that
\begin{align*}
	T_0 \ll& \paren{x^2+t^2}^{-1/4} e^{\frac{\pi}{2}t} x^{-1} \int_0^x \Min{1+\abs{\log y},y^{-1/2}} f_0\paren{t y/x} dy \ll \paren{x^2+t^2}^{-1/4} t^{-1} e^{\frac{\pi}{2}t} \log t.
\end{align*}

Note that for $0<z<2 t^{1/3}$, we have
\[ \abs{J_{it}(z)-\sum_{k=0}^N \frac{(-1)^k (z/2)^{it+2k}}{k! \Gamma(k+1+it)}} \le \sum_{k=N+1}^\infty \frac{t^{2k/3}}{k! \abs{\Gamma(1+k+it)}} \ll t^{-\frac{1}{2}-\frac{N+1}{3}} e^{\frac{\pi}{2}t}, \]
so we take $N=2$, giving
\begin{align}
\label{eq:BesselDervT2}
	T_{2a} \ll& x^{-1/2} t^{-3/2} e^{\frac{\pi}{2}t} \log t.
\end{align}

For the fifth term, we apply \eqref{eq:GlobalY0} and \eqref{eq:GlobalY1} so that
\begin{align*}
	T_{2b} \ll& x^{-1/2} e^{\frac{\pi}{2}t} \sum_{k=0}^N t^{-\frac{1}{2} -\frac{k}{3}} \sum_\pm \paren{t^{1/3} \abs{\int_0^1 e^{\pm i t^{1/3} y} y^{it+2k} f_0(y) dy}+\abs{\int_0^1 e^{\pm i t^{1/3} y} y^{it+2k} f_0'(y) dy}}.
\end{align*}
Repeated integration by parts implies
\[ \abs{\int_0^1 e^{\pm i t^{1/3} y} y^{it+2k} f_0(y) dy}+\abs{\int_0^1 e^{\pm i t^{1/3} y} y^{it+2k} f_0'(y) dy} \ll t^{-100}, \]
and we have
\begin{align}
\label{eq:BesselDervT3}
	T_{2b} \ll& x^{-1/2} t^{-3/2} e^{\frac{\pi}{2}t}.
\end{align}

For the second and third terms, we apply the asymptotic expansion of the $J$-Bessel function \eqref{eq:EMOTJitFirstTerm},
\begin{align*}
	e^{-\frac{\pi}{2} t} T_{1a}(C) \ll& \sum_{j=0}^{M_a} \abs{\int_0^x Y_0(x-y) \paren{y^2+t^2}^{-1/4-j/2} e^{i\omega(y,t)} f_1(y/C) \wtilde{a}_j(1+(y/t)^2) \frac{dy}{y}} \\
	& \qquad +t^{-1/2-\frac{M_a+1}{3}} x^{-1/2}, \\
	e^{-\frac{\pi}{2} t} T_{1b}(C) \ll& \sum_{j=0}^{M_b} \abs{\int_0^x Y_0(x-y) \paren{y^2+t^2}^{-1/4-j/2} e^{i\omega(y,t)} f_1((x-y)/C) \wtilde{a}_j(1+(y/t)^2) \frac{dy}{y}} \\
	& \qquad + \paren{x^2+t^2}^{-1/4} x^{-M_b-2} \paren{C^{1/2}+\abs{\log C}},
\end{align*}
using
\[ \int_{C/2}^{2C} \abs{Y_0(y)} dy \ll C^{1/2}+\abs{\log C}, \qquad \int_{C/2}^{2C} \abs{Y_0(x-y)} \frac{dy}{y^{k+1}} \ll x^{-1/2} C^{-k}. \]
So we may take $M_a=M_b=2$, but actually the trivial bound on $j \ge 1$ is sufficient for
\begin{align*}
	e^{-\frac{\pi}{2} t} T_{1a}(C) \ll& \abs{\int_0^x Y_0(x-y) \paren{y^2+t^2}^{-1/4} e^{i\omega(y,t)} f_1(y/C) \frac{dy}{y}} + t^{-1} \paren{x^2+t^2}^{-1/4} \log (x+t), \\
	e^{-\frac{\pi}{2} t} T_{1b}(C) \ll& \abs{\int_0^x Y_0(x-y) \paren{y^2+t^2}^{-1/4} e^{i\omega(y,t)} f_1((x-y)/C) \frac{dy}{y}} + t^{-1} \paren{x^2+t^2}^{-1/4} \log (x+t).
\end{align*}

Assume, for the moment that $C \ge t^\epsilon$ (this is only new for the $T_{1b}$ integral), then applying the asymptotic expansion \eqref{eq:LargeArgY0} of $Y_0(y)$,
\begin{align*}
	e^{-\frac{\pi}{2} t} T_{1a}(C) \ll& \sum_\pm \sum_{j = 0}^2 \abs{\int_0^x (x-y)^{-1/2-j} \paren{y^2+t^2}^{-1/4} e^{i\omega(y,t)\pm i y} f_1(y/C) \frac{dy}{y}} \\
	& \qquad + t^{-1} \paren{x^2+t^2}^{-1/4} \log (x+t), \\
	e^{-\frac{\pi}{2} t} T_{1b}(C) \ll& \sum_\pm \sum_{j = 0}^{\ceil{1/\epsilon}-1} \abs{\int_0^x (x-y)^{-1/2-j} \paren{y^2+t^2}^{-1/4} e^{i\omega(y,t)\pm i y} f_1((x-y)/C) \frac{dy}{y}} \\
	& \qquad + t^{-1} \paren{x^2+t^2}^{-1/4} \log (x+t).
\end{align*}
We have
\[ \frac{\partial}{\partial u} \omega(u,t) = -\frac{\sqrt{u^2+t^2}}{u}, \qquad \frac{\partial^j}{\partial u^j} \omega(u,t) \ll t^2 \frac{1+(t/u)^{j-2}}{u^{4-j} \paren{u^2+t^2}^{j-\frac{3}{2}}} \ll \frac{t^2}{u^2 (u+t)^{j-1}}, j \ge 2, \]
so by \cref{lem:BKY8.1} with $T=1$, $U=C$, $Q=C+t$, $Y= \frac{t^2(C+t)}{C^2}$, the integrals are $\BigO{t^{-1} \paren{x^2+t^2}^{-1/4}}$ (we have an $x^{-1/2}$ in both cases, just need to recover the $t^{3/2}$, hence just $t^\epsilon$) unless
\[ -\sqrt{1+(t/y)^2}\pm 1 \ll t^\epsilon \paren{\frac{t}{C\sqrt{C+t}}+C^{-1}}, \]
for some $y \in [C/2,3C/2]$.
This is impossible unless $t \ll C$, in which case we must have $\pm=+$ and $t^2/C \ll t^\epsilon \paren{t/\sqrt{C}+1}$ so that $x \gg C \gg t^{2-\epsilon}$.
The trivial bound in that case becomes
\[ T_{1a}(C), T_{1b}(C) \ll_\epsilon t^{\epsilon-1} \paren{x^2+t^2}^{-1/4} e^{\frac{\pi}{2} t}. \]

For $T_{1b}$, in case $x/t \le C < t^\epsilon$, we perform a single integration by parts on the integral
\begin{align*}
	& \int_0^x \frac{Y_0(x-y) f_1((x-y)/C)}{\paren{y^2+t^2}^{3/4}} \paren{-i \frac{\sqrt{y^2+t^2}}{y} e^{i\omega(y,t)}} dy \\
	&= \int_0^x \paren{\frac{3}{2}\frac{y Y_0(x-y) f\paren{\frac{x-y}{C}}}{y^2+t^2}+C^{-1} Y_0(x-y) f'\paren{\frac{x-y}{C}} - Y_1(x-y) f\paren{\frac{x-y}{C}}} \frac{e^{i\omega(y,t)} dy}{\paren{y^2+t^2}^{3/4}},
\end{align*}
and trivially bounding the result gives
\[ T_{1b}(C) \ll_\epsilon t^{\epsilon-1} \paren{x^2+t^2}^{-1/4} e^{\frac{\pi}{2} t} \]
in this case as well.

\section{An aside on the method of Bump, Friedberg and Goldfeld}
\label{sect:BFGmethod}
In their Proposition 8.1, BFG \cite{BFG} consider the behavior as $y_1,y_2$ become arbitrarily large in the Fourier coefficients of a Poincar\'e series with kernel function
\[ f_\text{BFG}(xyk) = \psi_m(x) \delta_{\substack{y_1 y_2 \le 1\\ y_1 \ge 1}} y_1^{2s_1+s_2} y_2^{2s_2+s_1} \exp\paren{-2\pi(m_1 y_1+m_2 y_1)}. \]
Provided $\Re(s_1)$ and $\Re(s_2)$ (hence also $\Re(2s_1+s_2)$ and $\Re(2s_2+s_1)$) are sufficiently large, the Poincar\'e series with this kernel and its Fourier coefficients will converge, but $f_\text{BFG}$ and hence the Poincar\'e series are not smooth, and hence there is no reason to believe the spectral expansion should apply.
Of course, this is easily fixed by replacing $\delta_{\substack{y_1 y_2 \le 1\\ y_1 \ge 1}}$ with a smooth approximation.

Since $f_\text{BFG}$ already has good decay at each $y_i=0,\infty$, one might consider replacing it with some function $f(xyk)=\psi_m(x)f(y)$ which is smooth and compactly supported in $y \in (\R^+)^2$ to bypass all concerns of convergence.
For such an $f$, we can apply the spherical case of the Whittaker expansion \eqref{eq:WhittExpand}, which reduces to
\begin{align*}
	f(g) =& \int_{\frak{a}^0_0} \what{f}(\mu) W^0(g,\mu) \sinmu^{0*}(\mu) d\mu, \\
	\what{f}(\mu) :=& \int_Y f(y) W^0(y,-\mu) dy.
\end{align*}
This allows us to apply the definition of the Bessel functions and Whittaker inversion in the proof of \cref{sect:FormalPoin}, so equality of the two interpretations of the Fourier coefficient looks like
\begin{equation}
\label{eq:SmoothedBFG}
\begin{aligned}
	& p_\rho(m) \int_{\mathcal{B}^{0*}} \frac{\lambda_\Xi(n) \wbar{\lambda_\Xi(m)}}{L(1,\AdSq \Xi)} \what{f}(\mu) W^0(y,\mu_\Xi) d_H\Xi \\
	& \qquad = \sum_{w \in \Weyl} \sum_{c \in \Z_{\ne 0}^2} S_w(m,n,c) p_{\rho^w}(n) F(y,mcwn^{-1} w^{-1}),
\end{aligned}
\end{equation}
where
\[ F(y,t) := \int_{\frak{a}^0_0} \what{f}(\mu) K^0_w(t,\mu) W^0(y,\mu) \sinmu^{0*}(\mu) d\mu. \]
Note that both the spectral side and every Weyl cell of the arithmetic side have a term $W^0(y,\mu)$, which decays exponentially as $y_1,y_2 \to \infty$ on compact subsets in $\mu$.
So every term in \eqref{eq:SmoothedBFG} is decaying to zero as $y_1,y_2 \to \infty$, and for $\Re(s_1),\Re(s_2)$ sufficiently large, this should also apply to the construction of BFG.
In fact, for smooth and compactly supported $f$, the function $\what{f}$ has super-polynomial decay in $\mu$, and hence every term of \eqref{eq:SmoothedBFG} has super-polynomial decay in $y_1, y_2$; we expect this to continue to hold for any smoothed (but not necessarily compactly supported) version of the BFG construction.

If we return to $f_\text{BFG}$ with the understanding that both the $w_5$ and $w_l$ terms are decaying in $y_1$ and $y_2$, we can find the difficulty by explicitly writing out the two terms.
We take $m=n=c=I$ for simplicity and note that the slight difference in the definitions of $w_l$ is not significant here.
Writing $f_y(g) = f_\text{BFG}(gy)$ and applying some simplifying substitutions (see the first paragraph of \cite[Sect 5]{Me01}), we have
\begin{align*}
	T_{w_l}(f_y)(I) =& y_1^{2-s_1-2s_2} y_2^{2-2s_1-s_2} \iiint\limits_{R_{w_l}(y)} \paren{1+x_1^2}^{-\frac{3s_2}{2}} \paren{1+x_2^2}^{-\frac{3s_1}{2}} \paren{1+x_3^2}^{\frac{1-3s_1-3s_2}{2}} \\
	& \qquad \times \exp-\paren{\frac{1}{y_2}\frac{\sqrt{1+x_1^2}}{\paren{1+x_2^2}\sqrt{1+x_3^2}}+\frac{1}{y_1}\frac{\sqrt{1+x_2^2}}{\paren{1+x_1^2}\sqrt{1+x_3^2}}} \\
	& \qquad \times \e{\frac{1}{y_1} \frac{x_1\sqrt{1+x_2^2}}{\paren{1+x_1^2}\sqrt{1+x_3^2}}+\frac{1}{y_2} \frac{x_1 x_3+x_2\sqrt{1+x_3^2}}{\paren{1+x_2^2}\sqrt{1+x_3^2}}} \\
	& \qquad \times \e{-y_1 \frac{x_2 x_3+x_1\sqrt{1+x_3^2}}{\sqrt{1+x_2^2}}-y_2 x_2} dx_1 \, dx_2 \, dx_3,
\end{align*}
\begin{align*}
	T_{w_5}(f_y)(I) =& y_1^{2-s_1-2s_2} y_2^{s_1-s_2} \iint\limits_{R_{w_5}(y)} \paren{1+x_1^2}^{\frac{1-3s_1-3s_2}{2}} \paren{1+x_3^2}^{-\frac{3s_2}{2}}\\
	& \qquad \times \exp-\paren{y_2 \frac{\sqrt{1+x_3^2}}{\sqrt{1+x_1^2}}+\frac{1}{y_1 y_2} \frac{1}{\sqrt{1+x_1^2}\paren{1+x_3^2}}} \\
	& \qquad \times \e{\frac{1}{y_1 y_2} \frac{x_3}{\sqrt{1+x_1^2}\paren{1+x_3^2}}+y_2\frac{x_1 x_3}{\sqrt{1+x_1^2}}-y_1 x_1} dx_1 \, dx_3,
\end{align*}
where
\begin{align*}
	R_{w_l}(y) =& \set{x \in \R^3 \setdiv \frac{1}{\paren{1+x_1^2} \sqrt{1+x_2^2} \paren{1+x_3^2}} \le y_1 y_2, \frac{\sqrt{1+x_1^2}}{\paren{1+x_2^2}\sqrt{1+x_3^2}} \ge y_2}, \\
	R_{w_5}(y) =& \set{(x_1,x_3) \in \R^2 \setdiv \frac{1}{\paren{1+x_1^2}\sqrt{1+x_3^2}} \le y_1, \frac{\sqrt{1+x_3^2}}{\sqrt{1+x_1^2}} \ge \frac{1}{y_2}}.
\end{align*}

The argument in \cite[Prop 8.1]{BFG} was simply that as $y_1,y_2 \to \infty$, $R_{w_l}(y,c) \to \phi$, while $R_{w_5}(y,c_1) \to \R^2$, and this is clearly true.
On the other hand, looking at the integrals themselves, the integrand of $T_{w_5}(f_y)$ clearly has exponential decay on any fixed compact subset of $\R^2$ as $y_2 \to \infty$.
We suggest that actually, for any smooth approximation of $f_\text{BFG}$, both $T_w(f_y)$ have super-polynomial decay and this is coming from the interplay of the exponential and oscillating factors (i.e. the exponential decay factor is very small at any stationary points of the oscillation).

\bibliographystyle{amsplain}

\bibliography{HigherWeight}

\end{document}